\documentclass[12pt,a4paper,reqno]{amsart}
\usepackage[english]{babel}
\usepackage{amssymb,latexsym,amsfonts,amsthm,upref,amsmath}
\usepackage[foot]{amsaddr}
\usepackage[margin=1in]{geometry} 
\usepackage[dvipsnames,x11names]{xcolor}
 \definecolor{myblue}{HTML}{003399}
\usepackage{hyperref}
\hypersetup{colorlinks,citecolor=myblue,filecolor=black,linkcolor=myblue,urlcolor=myblue}
\usepackage{enumerate}  
\usepackage{tikz}
\usepackage{float}
\usepackage{cleveref}
\usepackage{mathtools}
\usepackage{cite}
\usepackage{filecontents}
\makeatletter
\newcommand{\leqnomode}{\tagsleft@true}
\newcommand{\reqnomode}{\tagsleft@false}
\newcommand{\cev}[1]{\reflectbox{\ensuremath{\vec{\reflectbox{\ensuremath{#1}}}}}}
\makeatother
\newtheorem*{thm*}{Theorem}
\newtheorem*{lem*}{Lemma}
\newtheoremstyle{prim}{}{}{\normalfont}{}{\bfseries}{.}{ }{}
\newtheoremstyle{stil}{}{}{\slshape}{}{\bfseries}{.}{ }{}
\theoremstyle{stil}
\newtheorem{thm}{Theorem}[section]
\newtheoremstyle{defi}{}{}{}{}{\bfseries}{.}{ }{}
\theoremstyle{defi}
\newtheorem{defn}[thm]{Definition}
\theoremstyle{defi}
\newtheorem{rem}[thm]{Remark}
\theoremstyle{stil}
\newtheorem*{mthm*}{Main Theorem}
\newtheorem*{kor*}{Corollary}
\newtheorem{pro}[thm]{Proposition}
\theoremstyle{stil}
\newtheorem{lem}[thm]{Lemma}
\theoremstyle{stil}
\newtheorem{kor}[thm]{Corollary}
\theoremstyle{prim}
\newtheorem{ex}[thm]{Example}
\newenvironment{prf}{\noindent \textit{Proof.}}{\null\hfill$\qed$\hskip
2mm\vskip 2mm}

\newcommand{\Yd}{ {\rm Y}^+(\mathfrak{g}_{N})}
\newcommand{\Ydgl}{ {\rm Y}^+(\mathfrak{gl}_{N})}
\newcommand{\Y}{ {\rm Y}(\mathfrak{g}_{N})}
\newcommand{\Ylg}{ {\rm Y}(\mathfrak{gl}_{N})}
\newcommand{\X}{ {\rm X}(\mathfrak{g}_{N})}
\newcommand{\Xlg}{ {\rm X}(\mathfrak{gl}_{N})}
\newcommand{\Yggl}{ {\rm Y}(\mathfrak{gl}_{N})}
\newcommand{\Yg}{ \mathcal{Y}(\mathfrak{g}_{N})}
\newcommand{\Ygl}{ \mathcal{Y}(\mathfrak{gl}_{N})}
\newcommand{\YYg}{ \widetilde{\mathcal{Y}}(\mathfrak{g}_{N})}
\newcommand{\YYgl}{ \widetilde{\mathcal{Y}}(\mathfrak{gl}_{N})}
\newcommand{\BB}{  \mathcal{B}(N)}
\newcommand{\BBt}{  \widetilde{\mathcal{B}}(N)}
\newcommand{\CC}{\mathbb{C}}
\newcommand{\Sc}{\mathcal{S}}
\newcommand{\Dc}{\mathcal{D}}
\newcommand{\Lc}{\mathcal{L}}
\newcommand{\tm}{\mathbb{T}^+_{[n]}}
\newcommand{\tN}{\mathbb{T}^+_{[N]}}
\newcommand{\lN}{\mathbb{L}_{[N]}}
\newcommand{\bN}{\mathbb{B}_{[N]}}
\newcommand{\tn}{\mathbb{T}^+_{[m]}}
\newcommand{\tmn}{\mathbb{T}^+_{[n,m]}(u,v)}
\newcommand{\tjj}{\mathbb{T}^+_{[1,1]}(u,v)}
\newcommand{\Vc}{\mathcal{V}^c(\mathfrak{g}_N)}
\newcommand{\Vcgl}{\mathcal{V}^c(\mathfrak{gl}_N)}
\newcommand{\Vcc}{\mathcal{V}^{\mathrm{crit}}(\mathfrak{g}_N)}
\newcommand{\Vccgl}{\mathcal{V}^{\mathrm{crit}}(\mathfrak{gl}_N)}
\newcommand{\R}{\wvr{R}}

\newcommand{\RR}{\wnr{\wvr{R}}}
\newcommand{\RRR}{\widehat{R}}
\newcommand{\ZZ}{\mathbb{Z}}
\newcommand{\wvr}{\overline}
\newcommand{\wtld}{\widetilde}
\newcommand{\wnr}{\underline}
\newcommand{\vac}{\mathop{\mathrm{\boldsymbol{1}}}}

\newcommand{\tr}{\mathop{\mathrm{tr}}}

\newcommand{\ot}{\otimes}
\newcommand{\ts}{\hspace{1pt}}
\newcommand{\qdet}{ \mathop{\rm qdet} }
\newcommand{\sgn}{ \mathop{\rm sgn}}
\newcommand{\sdet}{ \mathop{\rm sdet}}

\newcommand{\ndo}{\mathop{\mathrm{End}}}
\newcommand{\om}{\mathop{\mathrm{Hom}}}

\newcommand{\Sym}{\mathfrak S}

\newcommand{\cdotrl}{\mathop{\hspace{-2pt}\underset{\text{RL}}{\cdot}\hspace{-2pt}}}
\newcommand{\cdotlr}{\mathop{\hspace{-2pt}\underset{\text{LR}}{\cdot}\hspace{-2pt}}}

\newcommand{\gl}{\mathfrak{gl}}
\newcommand{\g}{\mathfrak{g}}
\newcommand{\on}{\mathfrak{o}}
\newcommand{\spn}{\mathfrak{sp}}
\newcommand{\ve}{\varepsilon}
\newcommand{\kp}{\kappa}
\newcommand{\DY}{ {\rm DY}}
\newcommand{\Bc}{\mathcal{B}}

\newcommand{\fand}{\quad\text{and}\quad}
\newcommand{\Fand}{\qquad\text{and}\qquad}

\newcommand{\non}{\nonumber}
\newcommand{\beq}{\begin{equation}}
\newcommand{\eeq}{\end{equation}}
\newcommand{\ben}{\begin{equation*}}
\newcommand{\een}{\end{equation*}}

\makeatletter
\def\smalloverbrace#1{\mathop{\vbox{\m@th\ialign{##\crcr\noalign{\kern3\p@}%
  \tiny\downbracefill\crcr\noalign{\kern3\p@\nointerlineskip}%
  $\hfil\displaystyle{#1}\hfil$\crcr}}}\limits}
\makeatother

\makeatletter
\def\smallunderbrace#1{\mathop{\vtop{\m@th\ialign{##\crcr
   $\hfil\displaystyle{#1}\hfil$\crcr
   \noalign{\kern3\p@\nointerlineskip}%
   \tiny\upbracefill\crcr\noalign{\kern3\p@}}}}\limits}
\makeatother

\begin{document}

\title{Deformed quantum vertex algebra modules associated with braidings}

\author{Lucia Bagnoli}
\address[L. Bagnoli]{Dipartimento di Matematica, Sapienza Universit\`{a} di Roma, P.le Aldo Moro 5, 00185 Rome, Italy \& INFN sezione di Roma}
\email{lucia.bagnoli@uniroma1.it}

\author{Slaven Ko\v{z}i\'{c}}
\address[S. Ko\v{z}i\'{c}]{University of Zagreb Faculty of Science, Department of Mathematics, Bijeni\v{c}ka cesta 30, 10000 Zagreb, Croatia}
\email{slaven.kozic@math.hr}

\begin{abstract}
 We introduce the notion of deformed  quantum vertex algebra module    associated with   a braiding map. We construct two families of braiding maps over the Etingof--Kazhdan quantum vertex algebras  associated with the rational $R$-matrices of classical types. We investigate their properties and demonstrate the applications of the corresponding deformed modules to the   (generalized) Yangians and reflection algebras.
\end{abstract}

\maketitle

\allowdisplaybreaks

\section{Introduction}
\numberwithin{equation}{section}

A fundamental problem in the theory of vertex algebras, which goes back to I. Frenkel and Jing \cite{FJ}, is to associate suitable vertex algebra-like structures, {\em quantum vertex algebras} to different classes of quantum groups. It is motivated by the rich interplay between vertex algebras and representations of various infinite dimensional Lie algebras; see, e.g. the books by Kac \cite{Kac_book} and Lepowsky and Li \cite{LL_book}. 
Motivated by the ideas of E. Frenkel and  Reshetikhin \cite{FR},
Etingof and Kazhdan \cite{EK} introduced
the notion of  quantum vertex algebra. Furthermore, they   constructed first examples of these structures associated with the rational, trigonometric and elliptic $R$-matrix of type $A$. Later on, some other important generalizations of the notion of vertex algebra,  such as {\em field algebras} of Bakalov and Kac \cite{BKac}, {\em axiomatic $G_1$-vertex algebras}, inspired by Borcherds' {\em $G$-vertex algebras} \cite{Bor}, and {\em nonlocal vertex algebras} of Li \cite{LiG1,Li0}  and {\em $H_D$-quantum vertex algebras} of  Anguelova and Bergvelt  \cite{AB},   were  studied. A major step forward in the development of  quantum vertex algebra theory was Li's introduction of   {\em $\phi$-coordinated modules} \cite{Liphi} which, in   particular,   enabled the association of quantum vertex algebras with quantum affine algebras; see, e.g., the papers by   Jing,  Kong,   Li and Tan \cite{JKLT2} and Kong \cite{Kong}.
For a more detailed overview of the evolution of   quantum vertex algebra theory and its most   recent results see  \cite{BJK,BM,DGK,JKLT,Kong2} and references therein.

Suppose $(V, Y,\vac)$ is a vertex algebra and $W$ a $V$-module. The  corresponding module map $ Y_W(\cdot ,z)$ possesses two key properties (see, e.g., \cite{Kac_book,LL_book}), {\em locality}: for any $u,v\in V$
\beq\label{intro1}
Y_W(u,z_1)\ts Y_W(v,z_2) \sim
Y_W(v,z_2)\ts Y_W(u,z_1),
\eeq
 and the {\em weak associativity}:  for any $u,v\in V$ and $w\in W$
\beq\label{intro2}
Y_W(u,z_0+z_2)\ts Y_W(v,z_2)w \sim
Y_W(Y (u,z_0)v,z_2)w.
\eeq
The tilde ``$\sim$'' in \eqref{intro1} (resp. \eqref{intro2}) indicates that the given expressions coincide when multiplied by $(z_1-z_2)^r$ (resp. $(z_0+z_2)^r$) for  a sufficiently large positive integer $r$,   which depends on the choice of $u,v\in V$ (resp. $u,v\in V$, $w\in W$). The analogue of \eqref{intro1} for a quantum vertex algebra $(V, Y,\vac,\mathcal{S})$     and its module $W$ is  {\em $\mathcal{S}$-locality}: for any $u,v\in V$
\beq\label{intro3}
Y_W( z_1) \left( 1\ot Y_W( z_2) \right)\mathcal{S} (z_1-z_2)(u\ot v)\sim_h
Y_W(v, z_2) \ts Y_W(u, z_1)  ,
\eeq
where $\Sc$ is the {\em braiding of $V$}, a certain unitary map satisfying the quantum Yang--Baxter equation; cf. \cite{EK,Li}. On the other hand, the analogue of \eqref{intro2} for quantum vertex algebra modules, the
{\em weak associativity},   is of the same form:  for any $u,v\in V$ and $w\in W$
\beq\label{intro4}
Y_W(u,z_0+z_2)\ts Y_W(v,z_2)w \sim_h
Y_W(Y (u,z_0)v,z_2)w;
\eeq
cf. \cite{EK,Li}.
However, both $V$ and $W$ are  now     topologically free $\CC[[h]]$-modules and    the symbol ``$\sim_h$''
  indicates that  the form of $\mathcal{S}$-locality \eqref{intro3} (resp.    weak associativity \eqref{intro4})  generalizes its classical counterpart \eqref{intro1} (resp.    \eqref{intro2}) using the $h$-adic topology. The products of the given expressions with $(z_1-z_2)^r$ (resp. $(z_0+z_2)^r$) are required to coincide modulo $h^n$ only, and the power $r$ now depends on the choice of positive integer $n $ as well.

This paper is motivated in part by the notion of $H_D$-quantum vertex algebra    \cite{AB}, which    presents a generalization of  Etingof--Kazhdan's   quantum vertex algebras   \cite{EK}.
  In particular, it features a more general version  of $\mathcal{S}$-locality and a deformation of weak associativity  which is  
	   governed by a certain braiding map.  We study a new  concept of {\em $(\mu,\nu)$-deformed module}, or, more briefly, {\em deformed module} for a quantum vertex algebra $V$, such that its analogues of  \eqref{intro3} and \eqref{intro4} also depend on   certain  maps
\beq\label{munumaps}
\mu,\nu\colon V\otimes V\to V\otimes V\otimes\mathbb{C}[z^{\pm 1}][[x^{-1},h]].
\eeq
In addition, our notion of deformed module conforms to a slightly uncommon {\em truncation condition}: its module map  is allowed to possess essential singularities at zero while, on the other hand, it is required to be (topologically) truncated from above.  
The role of this   feature is to accommodate the well-known phenomenon which comes from the   theory of Yangians. Namely,
 the techniques and constructions from this framework can be often extended   to other classes of algebras associated with the same $R$-matrix, such as  generalized Yangians, twisted Yangians, reflection algebras etc.; see, e.g.,   \cite{KR,M,MNO,MR}. 
The truncation condition ensures that the theory of deformed modules  fits  all these algebras, so that some of their fundamental  properties naturally emanate from the    corresponding Etingof--Kazhdan quantum affine vertex algebra.
We  shall demonstrate this occurrence  in detail in the setting of   generalized  Yangians of Krylov and Rybnikov \cite{KR} and reflection algebras of Molev and Ragoucy \cite{MR}.

Recall the maps \eqref{munumaps}.		
In contrast with \eqref{intro3}, which depends on the original  braiding $\Sc$   of $V$, the deformed module map $Y_W(\cdot ,z)$  satisfies the {\em $\mu$-commutativity}: for any $u,v\in V$
\beq\label{intro5}
Y_W(z_1)\big(1\otimes Y_W(z_2)\big) \mu( -z_2+z_1,z_2)(u\otimes v) = Y_W(v,z_2)  Y_W(u,z_1)     .
\eeq
Roughly speaking, the map $\mu$  is   a unitary solution of the   quantum Yang--Baxter equation,  which is conjugated with  $\mathcal{S}$ by $\nu$, i.e. we have 
$$
\Sc(z)\ts \nu_{21}(-z,x+z)=\nu(z,x)\ts \mu(z,x).
$$
We refer to such a pair $(\mu,\nu)$ as a {\em compatible pair}.
As one may expect, such a structure does no  longer need to possess the weak associativity property \eqref{intro4}. However, we require     that it exhibits the {\em     $\nu$-associativity}:  for any $u,v\in V$ and $w\in W$
\beq\label{intro6}
 Y_W(u,z_2 +z_0)\ts Y_W(v,z_2)\ts w 
=    Y_W\big(Y ( z_0 )\ts \nu( z_0,z_2)(u\ot v),z_2\big)\ts w .
\eeq

Let us describe the organization and  main results of this manuscript.
In Section \ref{section_01}, we provide some preliminaries on the double Yangians  and quantum affine vertex algebras of classical types. In Section \ref{section_02}, we introduce and discuss the notion of  $(\mu,\nu)$-deformed module. Next, in Section \ref{section_03}, we construct two compatible pairs, $(\sigma,\rho)$ and    $(\mu,\nu)$.
The construction of $(\sigma,\rho)$ is carried out in the setting of Etingof--Kazhdan's quantum   vertex algebras $\Vc$ associated with rational $R$-matrices of types $A$--$D$ at an arbitrary level $c\in\CC$. Furthermore, we use the fusion procedure for these $R$-matrices to investigate   fixed points of $\sigma$ at the critical level $c=c_{\text{crit}}$. For $\mathfrak{g}_N$ of type $A$ and $B$, we show that the restriction of $\sigma$ on the tensor square of  {\em quantum Feigin--Frenkel center}, i.e. the center of $\Vcc=\mathcal{V}^{c_{\text{crit}}}(\mathfrak{g}_N) $, is the identity, while for types $C$ and $D$, we only obtain some partial results of the same form, due to  a lack of explicit formulae for complete sets of Segal--Sugawara vectors in these types; cf. \cite{JKMY,BJK}. 
Finally, the construction of  pair $(\mu,\nu)$ is given for the  quantum affine  vertex algebra $\Vccgl$ associated with rational $R$-matrix of type  $A$ at the critical level. Moreover, some families of fixed points of $\mu$, which come from  the quantum Feigin--Frenkel center, are found in this case as well.

As an application, in Section  \ref{section_06} (resp. Section \ref{section_07}),  we use the maps from Section \ref{section_03} to
 show that the  generalized Yangians (as well as   ordinary Yangians) of types $A$--$D$  (resp. certain reflection algebras) are naturally equipped with   structure of $(\sigma,\rho)$-deformed (resp. $(\mu,\nu)$-deformed) $\Vc$-module. Next, we establish an equivalence between certain wide classes of modules associated with such algebras and deformed $\Vc$-modules. Finally, we show that commutative subalgebras and families of central elements in these algebras naturally arise from the   families of fixed points of $\sigma$ and $\mu$   found in Section \ref{section_03}.

\section{Preliminaries}\label{section_01}
 
In this section, we recall some preliminary definitions and results on the double Yangians of types $A$--$D$ and the associated quantum vertex algebras.

\subsection{Double Yangians of classical types} 
Fix an integer $N\geqslant 2$. Let $\gl_N$ be the general linear, $\on_N$   the orthogonal  and   $\spn_N$   the symplectic Lie algebra, where $N$ is even in the symplectic case.
In order to consider all of them simultaneously, we denote by $\g_N$ any of these Lie algebras.
In the case  $\g_N=\on_N,\spn_N$,   introduce   the scalars $\ve_1=\ldots =\ve_N =1$ if $\g_N =\on_N$ and   $\ve_1=\ldots =\ve_{N/2} =1$, $\ve_{(N+2)/2}=\ldots =\ve_{N} =-1$ if $\g_N=\spn_N$. For any $i=1,\ldots ,N$ set $i^\prime=N-i+1$. Finally, let $G=(g_{ij})\in\ndo\CC^N$  be the matrix defined by $g_{ij}=\delta_{ij^\prime}\ve_i$ for all $i,j=1,\ldots ,N$. Clearly,
$G$ is symmetric in the  orthogonal and   skew-symmetric in the symplectic case.
For any matrix $A=(a_{ij})$ in $\ndo\CC^N$ let $A'=G A^t G^{-1}$, where $A^t$ denotes the transposed matrix $A^t=(a_{ji})$. We have $A^\prime=(\ve_i \ve_j a_{j'i'})$.
Set
$$
\ve =\begin{cases}
1,&\text{if }\g_N=\gl_N,\on_N,\\
2,&\text{if }\g_N=\spn_N.
\end{cases}
$$

Let $h$ be a formal parameter and $\CC[[h]]$ a commutative ring of formal power series in $h$.  Define the operators $I$, $P$ and $Q$ in $\ndo\CC^N \ot\ndo\CC^N$   by
$$I=\sum_{i,j=1}^N e_{ii}\ot e_{jj},\quad P=\sum_{i,j=1}^N e_{ij}\ot e_{ji}\fand Q=\sum_{i,j=1}^N \ve_{i}\ve_j \ts e_{ij}\ot e_{i^\prime j^\prime},$$
where $e_{ij}$ are matrix units. 
For $\g_N=\gl_N$, we consider the   {\em Yang $R$-matrix} $R(u)$ 
\beq\label{RA}
R(u)=I-\frac{h}{u}P
\eeq
and for $\g_N=\on_N,\spn_N$,   the {\em $R$-matrix}, which goes back to \cite{ZZ},
\beq\label{RBCD}
R(u)=I -\frac{h}{u}P + \frac{h}{u-h\kp}Q,
\quad\text{
where
}\quad
 \kp=\begin{cases}
N/2-1,&\text{if }\g_N=\on_N,\\
N/2+1,&\text{if }\g_N=\spn_N.
\end{cases}
\eeq

The $R$-matrices in \eqref{RA} and \eqref{RBCD} satisfy the {\em quantum Yang--Baxter equation}
\beq\label{RYBE}
R_{12}(u)\ts R_{13}(u+v)\ts R_{23}(v)=R_{23}(v)\ts R_{13}(u+v)\ts R_{12}(u).
\eeq
Both sides of \eqref{RYBE} are regarded as operators on the   tensor product $(\CC^N)^{\ot 3}$ and the subscripts
indicate the copies of $\CC^N$ on which the $R$-matrices act, e.g.,  $R_{12}(u)= R(u) \ot 1$. 
Let $\R(u)=f(u/h)R(u)$ be the {\em normalized $R$-matrix}, where for $\g_N=\gl_N$ (resp. $\g_N=\on_N,\spn_N$),    $f(u)$ is  a unique series in $1+u^{-1}\CC[[u^{-1}]]$ (resp. $1+u^{-2}\CC[[u^{-1}]]$) which satisfies
\begin{align*}
&f(u)\ts f(u+N)^{-1}= f(u)\ts f(-u)=\left(1-u^{-2}\right)^{-1}\\
\text{(resp. } &f(u)\ts f(u+\kp)=f(u)\ts f(-u)=\left(1-u^{-2}\right)^{-1}\text{ );}
\end{align*}
see   \cite[Subsect. 1.2]{BJK}  and \cite[Subsect. 2.2]{JKMY}   for more details on $f(u)$. The normalized $R$-matrices
 possess  the {\em unitarity}  property,
\beq\label{Runi}
\R(u)\ts \R(-u)=\R(-u)\ts \R(u) =1 ,
\eeq
and the {\em crossing symmetry}  property,\footnote{Note that we have
$\R(u)^{t_1}=\R(u)^{t_2}$ and $\R(u)^{\prime_1}=\R(u)^{\prime_2}$, so   these transposed $R$-matrices
 in \eqref{RCSYMA} and \eqref{RCSYMBCD} are denoted simply by $\R(u)^{t }$ and $\R(u)^{\prime }$, respectively.}
\begin{align}
&\R(-u)^t\ts \R(u+hN)^t =1\quad\text{ for }\g_N=\gl_N\label{RCSYMA}\\
&\R(u)\ts \R(u+h\kp  )^\prime =1\quad\text{ for }\g_N=\on_N,\spn_N.\label{RCSYMBCD}
\end{align}
Alternatively, the crossing symmetry property can be expressed as
\begin{align}
&\R(-u) \cdotlr \R(u+hN)  =\R(-u) \cdotrl \R(u+hN)  =1\quad\text{ for }\g_N=\gl_N,\label{RCSYMAx}\\
&\R(u)^\prime\cdotlr \R(u+h\kp  ) =\R(u)^\prime\cdotrl \R(u+h\kp  ) =1\quad\text{ for }\g_N=\on_N,\spn_N,\label{RCSYMBCDx}
\end{align}
where ``$\cdotlr$'' (resp. ``$\cdotrl$'') denotes the standard multiplication in $\ndo\CC^N\ot(\ndo\CC^N)^{\text{op}}$ (resp. $(\ndo\CC^N)^{\text{op}}\ot \ndo\CC^N $) and $A^{\text{op}}$  is the opposite algebra of $A$. In addition, we observe that for $\g_N=\on_N,\spn_N$, unitarity \eqref{Runi} and crossing symmetry \eqref{RCSYMBCDx} imply
$\R(-u)^{\prime}=\R(u+h\kp)$. This can be used to express the crossing symmetry alternatively as
\begin{align}
\R(-u+2h\kp) \cdotlr \R(u   ) =\R(-u+2h\kp) \cdotrl \R(u   ) =1\quad\text{ for }\g_N=\on_N,\spn_N.\label{RCSYMBCDxx}
\end{align}

Following \cite{I}, we define the  {\em double Yangian} $\DY(\gl_N)$      as the $h$-adically completed associative algebra over the ring $\CC[[h]]$ generated by the central element $C$ and the  elements $t_{ij}^{(\pm r)}$, where $i,j=1,\ldots ,N$ and $r=1,2,\ldots ,$ subject to the defining relations
\begin{align}
R(u-v)\ts T_1^{\pm}(u)\ts T_2^{\pm}(v)&=T_2^{\pm}(v)\ts T_1^{\pm}(u)\ts R(u-v),\label{DY1}\\
\R(u-v+h\varepsilon C/2)\ts T_1^-(u)\ts T_2^+(v)&=T_2^+(v)\ts T_1^-(u)\ts \R(u-v-h\varepsilon C/2).\label{DY2}
\end{align}
The generator matrices $T^{\pm}(u)$    are defined by
\beq\label{TDY}
T^-(u)=\sum_{i,j=1}^N e_{ij}\ot t^-_{ij}(u)\fand T^{+}(u)=\sum_{i,j=1}^N e_{ij}\ot t_{ij}^{+}(u) 
\eeq
and   the series   $t_{ij}^\pm (u) $ are given by
\beq\label{tDY}
t^-_{ij}(u)=\delta_{ij}+h\sum_{r=1}^{\infty} t_{ij}^{(r)}u^{-r}\fand t_{ij}^+ (u)=\delta_{ij}-h\sum_{r=1}^{\infty} t_{ij}^{(-r)}u^{r-1}.
\eeq
Next, following \cite{JLY}, we define the {\em double Yangian} $\DY(\g_N)$ for $\g_N=\on_N,\spn_N$   analogously, using \eqref{DY1}--\eqref{tDY}, but with two additional families of defining relations:
\beq\label{DY3}
 T^-(u) \ts T^-(u+h\kp)^\prime=1
\fand T^+(u)\ts  T^+(u+h\kp)^\prime=1.
\eeq

Let us recall two important subalgebras of the double Yangian. The {\em  Yangian} $\Y $ (resp. the {\em dual Yangian} $\Yd$) is defined as the $h$-adically completed subalgebra of $\DY(\g_N)$ generated by   all elements $t_{ij}^{(r)}$ (resp. $t_{ij}^{(-r)}$), where $i,j=1,\ldots ,N$ and $r=1,2,\ldots.$   By the Poincar\'{e}--Birkhoff--Witt theorem for the double Yangian   \cite[Thm. 2.2]{JKMY} and \cite[Thm. 3.4]{JLY},   $\Y $ (resp. $\Yd$) is isomorphic to the $h$-adically completed associative algebra over $\CC[[h]]$ generated by  the elements $t_{ij}^{(r)}$ (resp. $t_{ij}^{(-r)}$), where $i,j=1,\ldots ,N$ and $r=1,2,\ldots,$ subject to the defining relation  \eqref{DY1} for $T^-(u)$ (resp. for $T^+(u)$) and,   in the case $\g_N=\on_N,\spn_N$,
the first (resp. the second)   relation in \eqref{DY3}.

For any $c\in\CC$, define the {\em double Yangian $\DY_{\hspace{-1pt} c}(\g_N)$ at the level $c$} as the quotient of the algebra $\DY(\g_N)$ by its $h$-adically completed two-sided ideal generated by $C-c$. The {\em vacuum module $\Vc $ at the level  $c$} is defined as the $h$-adically  completed quotient
$
\DY_c(\g_N) / I
$, where $ I$ denotes the $h$-adically completed left ideal in $\DY_c(\g_N)$ generated by the elements $t_{ij}^{(r)}$ with $i,j=1,\ldots ,N$ and $r=1,2,\ldots .$ By the Poincar\'{e}--Birkhoff--Witt theorem for the double Yangian,  $\Vc $ and $\Yd$ are isomorphic as  $\CC[[h]]$-modules.

\subsection{Fusion procedure }\label{section_04x}

Let $R(u)$ be the Yang $R$-matrix  \eqref{RA}.
We shall recall a particular case of  the   fusion procedure   for $R(u)$   \cite{J};  see also \cite[Sect. 6.4]{M} for more details and references.
Recall that the symmetric group $\Sym_n$ acts on  $(\CC^N)^{\ot n}$ by permuting the tensor factors.
Denote by $ \mathcal{E}_{[n]}$   the action  of the  anti-symmetrizer $\frac{1}{n!}\sum_{p\in \Sym_n}\sgn p\cdot p \in \CC[\Sym_n]$   
on   $(\CC^N)^{\ot n}$. 
Consider the product 
\beq\label{r-product}
R_{[n]}(u ) =R_{[n]}(u_1,\ldots ,u_n) =\frac{1}{n!}\prod_{1\leqslant i<j\leqslant n}R_{ij}(u_i -u_j),
\eeq
where the factors are ordered lexicographically  on the pairs $(i,j)$.
Introduce the $n$-tuple
 \beq\label{uem0}
u_{[n]} =(u,u-h,\ldots ,u -(n-1) h),
\eeq
where $u$ is a single variable.
By \cite{J}, the consecutive evaluations $u_1 =0, u_2=-h, u_3=-2h,\ts \ldots ,\ts  u_n =-(n-1)h$ 
of the variables $u=(u_1,\ldots ,u_n)$ in 
\eqref{r-product}
 are well-defined and 
   they produce the following version of the {\em fusion procedure}:
\beq\label{fusionA}
R_{[n]}(u_1,\ldots ,u_n)\big|_{u_1=0}
\big|_{u_2=-h}\ldots \big|_{u_n=-(n-1)h}=  \mathcal{E}_{[n]}.
\eeq

	 Let $R(u)$ be the  $R$-matrix    \eqref{RBCD}.
Following the exposition in \cite{M}, we shall recall a special case of the fusion procedure for the Brauer algebra  \cite{B}; cf. \cite[Sect. 1.2]{Mnew}. 
Let $\omega$ be an indeterminate. Denote by $\Bc_n(\omega)$   the {\em Brauer algebra} over   $\CC(\omega)$ generated by   elements $s_1,\ldots ,s_{n-1}$ and $\varepsilon_1,\ldots ,\varepsilon_{n-1}$, which are subject to certain defining relations; see \cite[Sect. 3]{M} for details.
Its complex subalgebra  generated by  $s_1,\ldots ,s_{n-1}$ is isomorphic to the group algebra of   symmetric group $\Sym_n$ and the elements $s_i$ are identified with the transpositions $(i,i+1)$. 
Denote by $s_{ij}\in \Bc_n(\omega)$ with $i<j$  the element   corresponding to the transposition $(i,j)$ with respect to that isomorphism. 
Consider the elements  $\varepsilon_{ij}\in \Bc_n(\omega)$ defined by
$$\varepsilon_{j-1\ts j}=\varepsilon_{j-1}\Fand\varepsilon_{ij}= s_{i\ts j-1}\varepsilon_{j-1}s_{i\ts j-1}  \quad \text{for} \quad i<j-1 .$$
Let $s^{(n)}\in\Bc_n(\omega)$ be  the idempotent which corresponds to the one-dimensional representation of  $\Bc_n(\omega)$ which maps all $s_{ij}$ to the identity   and all $\varepsilon_{ij} $ to the zero. 
In the orthogonal case, denote by  $\mathcal{E}_{[n]}$ with $n=1,\ldots ,N$   the action of the idempotent $s^{(n)}\in \Bc_n (N)$ on  
$(\CC^N)^{\ot n}$   defined by
$$s_{ij}\mapsto P_{ij}\fand\varepsilon_{ij}\mapsto Q_{ij}\quad\text{for}\quad i<j.$$ In the symplectic case, denote by
$\mathcal{E}_{[n]}$ with $n=1,\ldots ,N/2$   the action of the idempotent $s^{(n)}\in\Bc_n (-N)$ on  
$(\CC^N)^{\ot n}$   defined by
$$s_{ij}\mapsto -P_{ij}\fand\varepsilon_{ij}\mapsto -Q_{ij}\quad\text{for}\quad i<j.$$
Let $R_{[n]}(u )=R_{[n]}(u_1,\ldots ,u_n) $ be the   $R$-matrix   product
of the same form as  
\eqref{r-product}.
For a single variable $u$, introduce the $n$-tuples
 \beq\label{uem1}
u_{[n]}=\begin{cases}
(u-(n-1)h,\ldots,u-h, u)&\text{for }\g_N=\on_N\text{ and }n=1,\ldots ,N,\\
(u,u-h,\ldots, u-(n-1)h)&\text{for }\g_N=\spn_N\text{ and }n=1,\ldots ,N/2.
\end{cases}
\eeq
By a special case of the  {\em fusion procedure}   for   $\Bc_n (\omega)$ \cite{IM,IMO},  we have
\beq\label{fusionBCD}
  R(u_{[n]}) = \mathcal{E}_{[n]}.
\eeq

\subsection{Etingof--Kazhdan's construction and the quantum Feigin--Frenkel center} 
Suppose $V$ is a $\CC[[h]]$-module. We shall denote by $V((z))_h$ (resp. $V[z^{\pm 1}]_h$) the $\CC[[h]]$-module of all power series $a(z)=\sum_{r\in\ZZ}a_r z^{-r}\in V[[z,z^{-1}]]$ (resp. $a(z)=\sum_{r\geqslant 0}a_r z^{\pm r}\in V[[z^{\pm 1} ]]$) such that $\lim_{r\to\infty}a_r =0$, where the limit is taken with respect to the $h$-adic topology.
Throughout the paper, we assume that the tensor products of $\CC[[h]]$-modules are  $h$-adically completed; see  \cite[Ch. XVI]{Kas} for more details on the $h$-adic topology.
For readers' convenience, we start  by recalling the definition of $h$-adic quantum vertex algebra;   see \cite[Sect. 1.4.1]{EK} and \cite[Def. 2.20]{Li}. In the definition, as well as in the rest of the paper, we use the usual vertex algebraic  expansion convention,   where the expressions such as
$(z_1+\ldots +z_n)^r$ with $r<0$ are expanded in nonnegative powers of $z_2,\ldots ,z_n$. Hence, in particular, we have $(z_1+z_2)^{r}\neq (z_2+z_1)^{r}$ for $r<0$.

\begin{defn}\label{qva}
Let $V$ be a topologically free $\mathbb{C}[[h]]$-module equipped with a distinguished element $\vac\in V$ (the {\em vacuum vector}), the $\CC[[h]]$-module map 
\begin{align}
Y(\cdot ,z) \colon V\ot V&\to V((z))_h,\label{Y}\\
u\ot v&\mapsto Y(z)(u\ot v)=Y(u,z)v=\sum_{r\in\mathbb{Z}} u_r v \ts z^{-r-1},\non
\end{align}
(the {\em vertex operator map})
and the $\CC[[h]]$-module map   
$$\Sc(z)\colon V\otimes V\to V\otimes V\otimes\mathbb{C}((z))[[h]]$$
(the {\em braiding}).
A    quadruple $(V,Y,\vac,\Sc)$ is said to be an {\em $h$-adic quantum vertex algebra} if it satisfies the following axioms.
\begin{enumerate}[(1)]
\item The vertex operator map possesses the   {\em weak associativity} property:
for any elements $u,v,w\in V$ and $n \in\mathbb{Z}_{\geqslant 0}$
there exists $r\in\mathbb{Z}_{\geqslant 0}$
such that
\begin{align*}
&\big((z_0 +z_2)^r\ts Y(u,z_0 +z_2)\ts Y(v,z_2)\ts w \big.\non\\
&\qquad\big. - (z_0 +z_2)^r\ts Y\big(Y(u,z_0)v,z_2\big)\ts w\big) 
\in  h^n   V[[z_0^{\pm 1},z_2^{\pm 1}]].
\end{align*}
\item The vacuum vector  satisfies the {\em vacuum property} and the {\em creation property}
$$
Y(\vac ,z)v=v,\quad Y(v,z)\vac \in V[[z]]\fand \lim_{z\to 0} Y(v,z)\ts\vac =v\quad\text{for all }v\in V.
$$
\item The  braiding  satisfies the {\em shift condition}
\begin{align}
&[\Dc\otimes 1, \mathcal{S}(z)]=-\frac{d}{dz}\mathcal{S}(z)\quad \text{for}\quad \Dc\in\ndo V\text{   defined by }\Dc v=v_{-2}\vac,\non\\
\intertext{the  {\em quantum Yang--Baxter equation}}
&\mathcal{S}_{12}(z_1)\ts\mathcal{S}_{13}(z_1+z_2)\ts\mathcal{S}_{23}(z_2)
=\mathcal{S}_{23}(z_2)\ts\mathcal{S}_{13}(z_1+z_2)\ts\mathcal{S}_{12}(z_1),\label{ybe}\\
\intertext{the {\em  unitarity condition}}
&\mathcal{S}_{21}(z)=\mathcal{S}^{-1}(-z),\label{uni}
\intertext{the {\em  hexagon identity}}
&\Sc(z_1)\left(Y(z_2)\ot 1\right) =\left(Y(z_2)\ot 1\right)\Sc_{23}(z_1)\ts \Sc_{13}(z_1+z_2),\label{hexagon}
\intertext{and the  {\em $\mathcal{S}$-locality}:
for any $u,v\in V$ and $n\in\mathbb{Z}_{\geqslant 0}$ there exists
$r\in\mathbb{Z}_{\geqslant 0}$ such that  }
 &\big((z_1-z_2)^{r}\ts Y(z_1)\big(1\otimes Y(z_2)\big)\big(\mathcal{S}(z_1 -z_2)(u\otimes v)\otimes w\big)\big. 
\non\\
 &\qquad\big. -(z_1-z_2)^{r}\ts Y(v,z_2) Y(u, z_1)  w \big)  \,
\in\,  h^n V[[z_1^{\pm 1},z_2^{\pm 1}]] \quad\text{for all }w\in V.\label{eslokaliti}
\end{align}
\end{enumerate}
\end{defn}

Let $m $ and $n$ be positive integers, $u=(u_1,\ldots ,u_n)$, $v=(v_1,\ldots, v_m)$ families of variables and $z$ a single variable.
Label the tensor factors of $(\ndo\mathbb{C}^{N})^{\ot  (n+m)}$
as follows,
\beq\label{nhk78}
\overbrace{(\ndo\mathbb{C}^{N})^{\ot  n}}^{1} \otimes
\overbrace{(\ndo\mathbb{C}^{N})^{\ot  m}}^{2}.
\eeq
Introduce the  formal power series   with coefficients in
\eqref{nhk78}
by
\begin{align}
\R_{nm}^{12}(u|v|z)= \prod_{i=1,\dots,n}^{\longrightarrow}
\prod_{j=n+1,\ldots,n+m}^{\longleftarrow} \R_{ij}(z+u_i -v_{j-n}),\label{rnm12}\\
\RR_{nm}^{12}(u|v|z)= \prod_{i=1,\dots,n}^{\longrightarrow}
\prod_{j=n+1,\ldots,n+m}^{\longrightarrow} \R_{ij}(z+u_i +v_{j-n}). \label{rnm123}
\end{align}
In \eqref{rnm12} and \eqref{rnm123}, the superscripts $1$ and $2$ indicate the tensor factors in \eqref{nhk78} while the arrows indicate the order of the factors. For example, we have
$$
\R_{22}^{12}(u|v|z)=\R_{14}\ts \R_{13} \ts \R_{24}\ts \R_{23}
\fand
\wnr{\wvr{R}}_{22}^{12}(u|v|z)= \RR_{13}\ts \RR_{14}\ts \RR_{23} \ts \RR_{24} 
$$
for $\R_{ij}=\R_{ij}(z+u_i -v_{j-n})$ and $\RR_{ij}=\R_{ij}(z+u_i +v_{j-n})$.
In order to simplify the notation, we write
$\R_{nm}^{12}(u|v)=\R_{nm}^{12}(u|v|0) $.
Furthermore, we shall use the notation 
\beq\label{teen}
T^\pm_{[n]}(u)=T_1^\pm(u_1)\ldots T_n^\pm(u_n )\fand T^\pm_{[n]} (u|z)=T^\pm_1(z+u_1)\ldots T^\pm_n(z+u_n),
\eeq
where the matrix entries of the coefficients of   \eqref{teen} are regarded as   operators on the vacuum  module $\Vc $. Hence, the coefficients of the expression  $T_{[n]}^\pm(u|z)$ belong  to 
$ (\ndo\CC^N)^{\ot n}\ot \ndo\Vc $.
Finally, we recall the Etingof--Kazhdan construction  \cite[Thm. 2.3]{EK}; see also \cite[Thm. 2.2]{BJK} and \cite[Thm. 4.1]{JKMY}. 
\begin{thm}\label{EK:qva}
For any $c\in \CC$
there exists a unique  structure of $h$-adic quantum vertex algebra
on  $\Vc$  such that the vacuum vector is
$\vac\in \Vc$, the vertex operator map $Y(\cdot ,z)$ is defined by
\beq\label{Y_form}
Y\big(T_{[n]}^+ (u)\vac,z\big)=T_{[n]}^+ (u|z)\ts T_{[n]}^- (u|z+h\varepsilon c/2)^{-1}
\eeq
  and the braiding $\mathcal{S} $ satisfies the identities  on
$
(\ndo\mathbb{C}^{N})^{\otimes n} \otimes
(\ndo\mathbb{C}^{N})^{\otimes m}\otimes \Vc \ot \Vc$, 
\begin{align}
&\mathcal{S}(z)\Big( T_{[m]}^{+24}(v) \ts
\R	_{nm}^{  12}(u|v|z-h\varepsilon   c) \ts T_{[n]}^{+13}(u)\ts \R_{nm}^{  12}(u|v|z)^{-1}\ts(\vac\otimes \vac) \Big) \non\\
=& \,\R_{nm}^{  12}(u|v|z)  \ts T_{[n]}^{+13}(u)  \ts\R_{nm}^{  12}(u|v|z+h\varepsilon   c)^{-1} \ts
 T_{[m]}^{+24}(v) \ts(\vac\otimes \vac) .\label{EK_S}
\end{align}
\end{thm}

Consider the $h$-adic quantum affine vertex algebras from Theorem \ref{EK:qva} at the {\em critical level} $c_{\text{crit}}$, where  $c_{\text{crit}}=-N$ for $\mathfrak{g}_N=\mathfrak{gl}_N$, while for $\mathfrak{g}_N=\mathfrak{o}_N$ (resp. $\mathfrak{g}_N=\mathfrak{sp}_N$) we have $c_{\text{crit}}=-N+2$ (resp. $c_{\text{crit}}=-\frac{N}{2}-1$). To simplify the notation, we shall denote $\mathcal{V}^{c_\text{crit}}(\mathfrak{g}_N)$ by $\Vcc$.
By \cite[Thm. 2.4]{BJK} and \cite[Thm. 4.9]{JKMY}, all coefficients of the series
\beq\label{tm}
\tm (u) =\tr_{1,\ldots,n} \mathcal{E}_{[n]}\ts T_{[n]}^+(u_{[n]}) \vac \in \Vcc[[u]],
\eeq
where $n=1,\ldots ,N$ (resp. $n=1,\ldots ,N/2$) for $\mathfrak{g}_N=\mathfrak{gl}_N,\mathfrak{o}_N$ (resp. $\mathfrak{g}_N=\mathfrak{sp}_N$),
belong to the {\em center} 
$$\mathfrak{z}(\Vcc)
=\left\{w\in \Vcc\,:\, Y(v,z)w\in \Vcc[[z]]\text{ for all }v\in\Vcc\right\}
$$ 
of   $\Vcc$; see \cite[Subsect. 3.2]{JKMY} for more information on the notion of center of $h$-adic quantum vertex algebra. 
Furthermore, by \cite[Thm. 2.8]{BJK} and
\cite[Thm. 4.9]{JKMY} we have 
	
	\begin{thm}\label{thm_ff}
	The center $\mathfrak{z}(\Vcc)$ for $\mathfrak{g}_N=\mathfrak{gl}_N$ and 
  for $\mathfrak{g}_N=\mathfrak{o}_N$ with $N$ odd (i.e. for $\mathfrak{g}_N$ of type  $A$ and $B$) 
	is a commutative associative algebra with respect to the multiplication $a\cdot b=a_{-1}b$. It is topologically generated by the coefficients of $\tm (u)$, $n=1,\ldots ,N$.
	\end{thm}
	
We remark that for $\mathfrak{g}_N=\mathfrak{o}_N,\spn_N$ with $N$ even (i.e. for $\mathfrak{g}_N$ of type  $C,D$) the coefficients of $\tm (u)$, $n=1,\ldots ,N/2$ again topologically generate a commutative associative algebra with respect to the multiplication $a\cdot b=a_{-1}b$. However, in this case, such an algebra does not exhaust the entire center $\mathfrak{z}(\Vcc)$; see \cite[Cor 2.10]{BJK}.

\section{Deformed modules for \texorpdfstring{$h$}{h}-adic quantum vertex algebras}\label{section_02}

The goal of this section is to present the notion of deformed module.
The next definition is
motivated by the properties \eqref{ybe} and \eqref{uni} of the braiding. Although it envolves an $h$-adic quantum vertex algebra, it  does not  depend on its structure,  so, alternatively, it can be given  over any   $\CC[[h]]$-module.

\begin{defn}\label{def27}
Let $(V,Y,\vac,\Sc)$  be an $h$-adic quantum vertex algebra. The map
$$\mu(z,x)\colon V\otimes V\to V\otimes V\otimes\mathbb{C}[z^{\pm 1}][[x^{-1},h]]$$
is said to be a {\em braiding map on $V$}, if it satisfies {the {\em quantum Yang--Baxter  equation}}
\begin{align}
&\mu_{12}(z_1,x+z_2)\ts\mu_{13}(z_1+z_2,x)\ts\mu_{23}(z_2,x)
=\mu_{23}(z_2,x)\ts\mu_{13}(z_1+z_2,x)\ts\mu_{12}(z_1,x+z_2), \label{gen_br_ybe}\\
\intertext{    the {\em  unitarity condition}}
&\mu_{21}(z,x-z)=\mu^{-1}(-z,x)\label{gen_br_uni}
\intertext{ and the {\em  singularity constraint}}
&\mu (-x+z,x) (a\ot b)\in V\ot V \otimes\mathbb{C}[z ][[x^{-1},h]] \text{ for all }a,b\in V.\label{gen_br_sing}
\end{align}
\end{defn}

\begin{rem}
Roughly speaking,   the singularity constraint \eqref{gen_br_sing} ensures that the map $\mu =\mu(z,x)$ does not possess a singularity at $z=-x$. For example, the series
$$
m(z,x)=\frac{1}{x+z}=\sum_{k\geqslant 0} (-1)^k\frac{z^k}{x^{k+1}}\in\CC[z][[x^{-1}]]\subset \mathbb{C}[z^{\pm 1}][[x^{-1},h]]
$$
does not satisfy such a requirement as $m(-x+z,x)=z^{-1}\not\in\mathbb{C}[z ][[x^{-1}]]$. We explain the role of this property   after Definition \ref{defn_deformed_top} below.
\end{rem}

Let $(V,Y,\vac,\Sc)$  be an $h$-adic quantum vertex algebra and
\beq\label{remark_on_nu}
\nu(z,x)\colon V\otimes V\to V\otimes V\otimes\mathbb{C}[z^{\pm 1}][[x^{-1},h]]
\eeq
an arbitrary map.
By  \eqref{Y} and \eqref{remark_on_nu}, the composition  $Y(z)\nu(z,x)$ is   well-defined and its image belongs to 
$V((z))[[x^{-1},h]]$.
To simplify the notation, we shall denote this composition  by
$Y^\nu(z,x)=Y^\nu(\cdot ,z,x)$. 
In what follows, we want to consider the braiding maps over  $h$-adic quantum vertex algebras
  which differ from the original braiding $\Sc$. They  will need to satisfy certain compatibility conditions, related to the structure of $V$, as follows.

\begin{defn}\label{def29}
Let $(V,Y,\vac,\Sc)$  be an $h$-adic quantum vertex algebra. A   braiding map $\mu=\mu(z,x)$
 on $V$ is said to be {\em compatible} if there exists a map  (the {\em   intertwining map})
$\nu(z,x)$ as in \eqref{remark_on_nu}
which satisfies
\begin{align}
& \Sc(z)\ts \nu_{21}(-z,x+z)=\nu(z,x)\ts \mu(z,x),\label{mu1}\\
& Y^\nu(a,z,x)b \in V[x^{-1}][[z^{\pm 1},h]]\text{ for all }a,b\in V.\label{mu2}
\end{align}
Furthermore, the pair $(\mu, \nu)$ is then said to be a  {\em  compatible pair}. 
\end{defn}

We discuss the role of the   requirement \eqref{mu2} after Definition \ref{defn_deformed_top} below.

\begin{ex}\label{trivial}
One easily deduces from   Definition  \ref{def29}  that the original braiding $\Sc$ of the $h$-adic quantum vertex algebra $(V,Y,\vac,\Sc)$ is compatible. Indeed, the corresponding intertwining map  is the embedding
$\iota\colon v_1\ot v_2\mapsto v_1\ot v_2\ot 1$. Clearly, the  both maps of the compatible pair $(\Sc,\iota)$ are constant with respect to the second variable $x$.
\end{ex}

 Our next goal is to introduce $(\mu,\nu)$-deformed modules. In comparison with the notion of $h$-adic quantum vertex algebra, they   employ  a slightly  more general topological setting.
Let $W$ be a $\CC[[h]]$-module such that there exists a family of $\CC[[h]]$-submodules
\beq\label{top_submodules}
W=W_0\supseteq W_1\supseteq W_2\supseteq \ldots
\qquad\text{
satisfying}
\qquad
h^n W\subseteq W_n\quad\text{for all }n\geqslant 0.
\eeq
We shall consider the linear topology $\tau$ on $W$ defined by the basis
\beq\label{top_basis}
w+W_n\quad\text{for all }n\geqslant 0,\, w\in W.
\eeq
Denote by $W((z^{-1}))_\tau$ the $\CC[[h]]$-module of all formal power series 
$$a(z)=\sum_{r\in\ZZ}a_r z^r\in W[[z^{\pm 1}]]\quad\text{such that}\quad  \lim_{r\to \infty} a_r=0,$$
where the limit is taken with respect to the topology $\tau$.  Such a notation naturally extends to the multiple variable case. For example, $a(z_1,z_2)\in W((z_1^{-1},z_2^{-1}))_\tau$ means that  the series $a=a(z_1,z_2)$ belongs to $ W[[z_1^{\pm 1},z_2^{\pm 1}]]$ and that for every  $n \in\mathbb{Z}_{\geqslant 0}$ it decomposes as
$a=b+c$
for some 
$b=b(z_1,z_2)\in W((z_1^{-1},z_2^{-1}))$ and  
$c=c(z_1,z_2)\in W_n[[z_1,z_2]]$.

\begin{defn}\label{defn_deformed_top}
Let $(V,Y,\vac,\Sc)$  be an $h$-adic quantum vertex algebra equipped with a
  compatible pair $(\mu,\nu)$.
 Let $W$ be a $\CC[[h]]$-module equipped with a descending family of submodules as in \eqref{top_submodules} such that it is complete with respect to the topology $\tau$ 
defined by the basis  \eqref{top_basis}. 
Let $Y_W(\cdot, z)$ be a $\CC[[h]]$-module map
\begin{align*}
Y_W(\cdot ,z) \colon V\ot W&\to W((z^{-1}))_\tau, \\
v\ot w&\mapsto Y_W(z)(v\ot w)=Y_W(v,z)w=\sum_{r\in\mathbb{Z}} v_r w \ts z^{-r-1}.\non
\end{align*}
A pair $(W,Y_W)$ is said to be a {\em   $(\mu,\nu)$-deformed $V$-module} if the   map $Y_W(\cdot, z)$ satisfies the {\em vacuum property}
\begin{align}
&Y_W(\vac, z)w=w\quad\text{for all }w\in W,\label{vac_cond_deformed}\\
\intertext{the {\em truncation condition}: for any $u,v\in V$ and $w\in W$   we have}
&Y_W(u,z_1)\ts Y_W(v,z_2)w\in W((z_1^{-1},z_2^{-1}))_\tau,\label{trunc_cond}
\intertext{the   {\em $\nu$-associativity}:
for any elements $u,v\in V$ and $w\in W$ 
we have}
 &   Y_W(u,z_2 +z_0)\ts Y_W(v,z_2)\ts w 
=    Y_W\big(Y^\nu(u,z_0,z_2)v,z_2\big)\ts w  ,\label{rho_assoc_top}
\intertext{and the {\em $\mu$-commutativity}: for any $u,v\in V$ and $w\in W$   we have  }
&   Y_W(z_1)\big(1\otimes Y_W(z_2)\big)\big(\mu( -z_2+z_1,z_2)(u\otimes v)\otimes w\big)= Y_W(v,z_2)  Y_W(u,z_1)  w     
.\label{sigma_loc_top}
\end{align}
 \end{defn}

Regarding Definition \ref{defn_deformed_top}, the expression on the left-hand side of $\nu$-associativity \eqref{rho_assoc_top} is well-defined due to  the truncation condition \eqref{trunc_cond}, while the right-hand side is well-defined due to \eqref{mu2}.
Furthermore, the singularity constraint \eqref{gen_br_sing} and the truncation condition  \eqref{trunc_cond}    ensure  that the left-hand side of   $\mu$-commutativity \eqref{sigma_loc_top} is well-defined.

\begin{defn}\label{def_gamma}
Let $(V,Y,\vac,\Sc)$  be an $h$-adic quantum vertex algebra equipped with  
  compatible pairs $(\mu,\nu)$ and $(\sigma,\rho)$. The pairs  $(\mu,\nu)$ and $(\sigma,\rho)$ are said to be {\em equivalent} if there
	exist    maps
$$ \gamma^{\pm 1}(z,x)\colon V\otimes V\to V\otimes V\otimes\mathbb{C}[z, z^{-1}][[x^{-1},h]],$$
such that
$
\gamma(z,x) \ts \gamma^{-1}(z,x) =\gamma^{-1}(z,x)\ts   \gamma(z,x)=1
$, 
which satisfy the {\em shift condition}
\begin{align}
&\gamma^{\pm 1}(z,x)=\gamma^{\pm 1}(-z,x+z) \label{shift_gamma}
\end{align}
and the equalities
\beq\label{main_gamma}
\mu(z,x)=\gamma^{-1}(z,x)\ts \sigma(z,x)\ts \gamma_{21}(z,x)\Fand
\nu(z,x)=\rho(z,x)\ts \gamma(z,x)
.
\eeq
	\end{defn}

\begin{rem}\label{rem_useful} 
The next  two simple observations will be useful later on.
The unitarity condition  \eqref{gen_br_uni} for the map $\sigma$ and the shift condition \eqref{shift_gamma} for the map $\gamma$ imply that the map $\mu$, defined by the first identity  in \eqref{main_gamma}, also satisfies the unitarity  \eqref{gen_br_uni}. 
Moreover, the property \eqref{mu1} for the pair $(\sigma,\rho)$ 
and the shift condition \eqref{shift_gamma} for the map $\gamma$
imply that the map $\nu$, defined by the second identity  in
   \eqref{main_gamma}, also possesses the property \eqref{mu1}.
\end{rem}

\section{Compatible braiding maps over \texorpdfstring{$h$}{h}-adic quantum affine vertex algebras}\label{section_03}

In this section, we construct two compatible pairs over the $h$-adic quantum vertex algebras from Theorem \ref{EK:qva}  and investigate their properties.

\subsection{Compatible pair \texorpdfstring{$(\sigma,\rho)$}{(sigma,rho)} over \texorpdfstring{$\Vc$}{Vc(gN)}}

The goal of this subsection is to construct an example of a compatible braiding map on $\Vc$ which differs from \eqref{EK_S}. The candidate is the map $\sigma=\sigma(z,x)$ from the following lemma, which is constant with respect to the second variable $x$; cf. Definition \ref{def27}.

\begin{lem}\label{lem_sigma}
There exists a unique map 
$$
 \sigma(z)\colon \Vc\ot\Vc\to\Vc\ot\Vc\ot\CC[[h]][z^{-1}]_h
$$
 such that on
$
(\ndo\mathbb{C}^{N})^{\otimes n} \otimes
(\ndo\mathbb{C}^{N})^{\otimes m}\otimes \Vc \ot \Vc$
we have
\beq\label{sigma_formula}
 \sigma(z)\Big( T_{[n]}^{+13}(u) \ts T_{[m]}^{+24}(v) (\vac\ot\vac)\Big)  
= \R_{nm}^{  12}(u|v|z)\ts  T_{[n]}^{+13}(u)  \ts T_{[m]}^{+24}(v) 
\R_{nm}^{  12}(u|v|z)^{-1}(\vac\ot\vac).
\eeq
Moreover, the map $\sigma$ is a braiding map on  $\Vc$.
\end{lem}

\begin{prf}
Lemmas \ref{lem_sigma}, \ref{lem_rho} and Proposition \ref{thm_compatiblity} can be verified by straightforward calculations. They are carried out in parallel to the   proofs of \cite[Thm. 2.2]{BJK}, \cite[Thm. 2.3.8]{G} and \cite[Thm. 4.1]{JKMY} and they employ   properties \eqref{RYBE}--\eqref{RCSYMBCD} of the $R$-matrix   and the   defining relations \eqref{DY1}, \eqref{DY2} and \eqref{DY3}. Therefore, we omit the proofs of the lemmas and only   sketch the proof of Proposition \ref{thm_compatiblity}.
\end{prf}

In order to show that $\sigma$ is compatible, we need to find its intertwining map.

\begin{lem}\label{lem_rho}
There exists a unique map 
$$
 \rho(z)\colon \Vc\ot\Vc\to\Vc\ot\Vc\ot\CC [[h]][z^{-1}]_h
$$
 such that on
$
(\ndo\mathbb{C}^{N})^{\otimes n} \otimes
(\ndo\mathbb{C}^{N})^{\otimes m}\otimes \Vc \ot \Vc$
we have
\beq\label{rho_formula}
 \rho(z)\Big( T_{[n]}^{+13}(u) \ts T_{[m]}^{+24}(v)(\vac\ot\vac) \Big)  
=   T_{[n]}^{+13}(u)  \ts \R_{nm}^{  12}(u|v|z+h\varepsilon c)^{-1}\ts T_{[m]}^{+24}(v) 
\R_{nm}^{  12}(u|v|z) (\vac\ot\vac).
\eeq
\end{lem}

Finally, it remains to show that the maps $\sigma$ and $\rho$ satisfy  Definition \ref{def29}.

\begin{pro}\label{thm_compatiblity}
The maps  $\sigma$ and  $\rho$ form a compatible pair on $\Vc$.
\end{pro}

\begin{prf}
In order to verify \eqref{mu1}, it is useful to  observe that the explicit formula for the action of $\rho_{43}(-z)$
on $
(\ndo\mathbb{C}^{N})^{\otimes n} \otimes
(\ndo\mathbb{C}^{N})^{\otimes m}\otimes \Vc \ot \Vc$
 is given by
$$
 \rho_{43}(-z)\Big( T_{[n]}^{+13}(u) \ts T_{[m]}^{+24}(v) (\vac\ot\vac)\Big)  
=   T_{[m]}^{+24}(v)  \ts \R_{nm}^{  12}(u|v|z-h\varepsilon c) \ts T_{[n]}^{+13}(u) 
\R_{nm}^{  12}(u|v|z)^{-1}(\vac\ot\vac) .
$$
By applying the braiding $\Sc (z) $ on the above equality, one finds that the image of
the expression
$T_{[n]}^{+13}(u)   T_{[m]}^{+24}(v)(\vac\ot\vac) $ under the composition
$\Sc (z)\rho_{43}(-z)$  
coincides with
 the right-hand side of \eqref{EK_S}. On the other hand, by using formulae \eqref{sigma_formula} and \eqref{rho_formula}, one easily checks that  the image of
the same expression  under the composition $\rho(z)\sigma(z)$ coincides with
 the right-hand side of \eqref{EK_S} as well. Thus, we conclude that  \eqref{mu1} holds.
The remaining property \eqref{mu2} obviously holds as the map $\rho=\rho(z,x)$ is constant with respect to the second variable $x$. 
\end{prf}

The maps   $\sigma$ and $\rho$ satisfy the following {\em hexagon relations}; cf. \eqref{hexagon}.

\begin{pro}
We have
\begin{align}
& \rho(z_1)\left(Y^\rho(z_2)\ot 1\right)=\left(Y^\rho(z_2)\ot 1\right) \rho_{13}(z_1+z_2)\ts \rho_{23}(z_1),\label{hex1}\\
& \sigma(z_1)\left(Y^\rho(z_2)\ot 1\right)=\left(Y^\rho(z_2)\ot 1\right)   \sigma_{23}(z_1)\ts \sigma_{13}(z_1+z_2).\label{hex2}
\end{align}
\end{pro}

\begin{prf}
Let us prove     \eqref{hex1}.
First, we observe that the map $Y^\rho$ satisfies the identity
\beq\label{Yrhoo}
Y^\rho(T_{[n]}^{+13 }(u),z_2)T_{[m]}^{+24 }(v)  (\vac\ot\vac ) 
=T_{[n]}^{+13 }(u|z_2) \ts T_{[m]}^{+23 }(v)  \vac .
\eeq
Indeed, this is easily verified by the use of   defining relation \eqref{DY2} for the double Yangian and formula \eqref{Y_form} for the vertex operator map.
 Clearly,  \eqref{Yrhoo} implies
\beq\label{Yrho}
Y^\rho(T_{[n]}^{+14 }(u),z_2)T_{[m]}^{+25 }(v)\ts T_{[k]}^{+36}(w) (\vac\ot\vac\ot\vac) 
=T_{[n]}^{+14 }(u|z_2) \ts T_{[m]}^{+24 }(v)\ts T_{[k]}^{+35}(w)(\vac\ot\vac ).
\eeq
 By applying $\rho(z_1)$ to \eqref{Yrho} and then using \eqref{rho_formula}, we find that  $\rho(z_1)\left(Y^\rho(z_2)\ot 1\right)$ maps
\beq\label{the_expression}
T_{[n]}^{+14 }(u)\ts T_{[m]}^{+25 }(v)\ts T_{[k]}^{+36}(w)(\vac\ot\vac\ot\vac)
\eeq
to
\begin{align}
&T_{[n]}^{+14}(u|z_2)  \ts 
T_{[m]}^{+24}(v) 
\R_{mk}^{  23}(v|w|z_1+h\varepsilon c)^{-1}\ts
\R_{nk}^{  13}(u|w|z_1+z_2+h\varepsilon c)^{-1}\non\\ 
&\qquad\times T_{[k]}^{+35}(w)\ts
\R_{nk}^{  13}(u|w|z_1+z_2 ) \ts 
\R_{mk}^{  23}(v|w|z_1 )(\vac\ot\vac ), \label{proof_prop34} 
\end{align}
where the superscripts $1$--$6$ in \eqref{the_expression} indicate the tensor factors as follows:
\beq\label{tfactors_sup}
\overbrace{(\ndo\mathbb{C}^{N})^{\ot  n}}^{1} \otimes
\overbrace{(\ndo\mathbb{C}^{N})^{\ot  m}}^{2}\otimes
\overbrace{(\ndo\mathbb{C}^{N})^{\ot  k}}^{3}\otimes
\overbrace{\Vc}^{4}
\otimes
\overbrace{\Vc}^{5}\otimes \overbrace{\Vc}^{6}.
\eeq
Finally, by using  \eqref{rho_formula} and then \eqref{Yrho},  one shows  that the image of \eqref{the_expression}
under the composition
$\left(Y^\rho(z_2)\ot 1\right) \rho_{13}(z_1+z_2)\ts \rho_{23}(z_1)$
equals \eqref{proof_prop34} as well, thus proving the hexagon relation   \eqref{hex1}. 
The remaining  relation \eqref{hex2} for $\sigma$    can be verified by   analogous arguments which rely on  \eqref{sigma_formula} and \eqref{Yrho}.
\end{prf}

\subsection{Compatible pair \texorpdfstring{$(\mu,\nu)$}{(mu,nu)} over \texorpdfstring{$\Vccgl$}{Vcrit(glN)}}

 The goal of this subsection is to construct another
compatible pair $(\mu,\nu)$  on $\Vccgl$, which is equivalent to $(\sigma, \rho)$.   We start by  fixing the auxiliary map $\gamma$ from   Definition \ref{def_gamma}.  

\begin{lem}\label{lem_gamma}
There exists a unique map 
$$
 \gamma(z,x)\colon \Vccgl\ot\Vccgl\to\Vccgl\ot\Vccgl\ot\CC[z^{\pm 1}][[x^{-1},h]]
$$
 such that on
$
(\ndo\mathbb{C}^{N})^{\otimes n} \otimes
(\ndo\mathbb{C}^{N})^{\otimes m}\otimes \Vccgl \ot \Vccgl$
we have
\begin{align}
 \gamma(z,x)\Big( T_{[n]}^{+13}(u) \ts T_{[m]}^{+24}(v) (\vac\ot\vac)\Big)  
=     T_{[n]}^{+13}(u)\ts  \RR_{nm}^{  12}(u|v|2x+z-hN)^{-1}  \ts T_{[m]}^{+24}(v)  (\vac\ot\vac).\label{gamma_form}
\end{align}
Its inverse is given by
\begin{align}
 \gamma^{-1}(z,x)\Big( T_{[n]}^{+13}(u) \ts T_{[m]}^{+24}(v) (\vac\ot\vac)\Big)  
=     T_{[n]}^{+13}(u) \ts \RR_{nm}^{  12}(u|v|2x+z ) \ts   T_{[m]}^{+24}(v)  (\vac\ot\vac) \label{gamma-_form}
\end{align}
and both maps satisfy the shift condition \eqref{shift_gamma}.
\end{lem}

\begin{prf}
One can show that the maps $\gamma$ and $\gamma^{-1}$ are well-defined  by arguing as in \cite[Lemma 2.1]{EK}.
The expression for the inverse \eqref{gamma-_form} is easily verified by using the crossing symmetry \eqref{RCSYMAx}.
Finally, the last assertion of the lemma is evident from \eqref{gamma_form} and \eqref{gamma-_form}.
\end{prf}

Motivated by  the first identity in \eqref{main_gamma}, we define the map $\mu$ by
\beq\label{mu_form_second}
\mu(z,x)=\gamma^{-1}(z,x)\ts \sigma(z)\ts \gamma_{21}(z,x).
\eeq
 Lemmas \ref{lem_sigma} and \ref{lem_gamma} imply that its image belongs to
$ \Vccgl^{\ot 2}\ot\CC[z^{\pm 1}][[x^{-1},h]]$. In fact, one can derive from these lemmas the explicit expression for the action of $\mu$,
\begin{align}
 \mu(z,x)\Big( T_{[n]}^{+13}(u) \ts T_{[m]}^{+24}(v)(\vac\ot\vac) \Big)  
= \,&\RR_{nm}^{  12}(u|v|2x+z-hN)^{-1}\cdotlr\left(\R_{nm}^{  12}(u|v|z)\ts  T_{[n]}^{+13}(u) \non\right.\\ 
  &\hspace{-25pt}\Big.\times\RR_{nm}^{  12}(u|v|2x+z )\ts  T_{[m]}^{+24}(v) \ts
\R_{nm}^{  12}(u|v|z)^{-1}  \Big)(\vac\ot\vac),\label{mu_form}
\end{align}
where ``$\cdotlr$'' denotes
the standard multiplication in 
$ (\ndo\CC^N)^{\ot n} \ot ((\ndo\CC^N)^{\text{op}})^{\ot m}$.

\begin{lem}\label{lem_mu}
The map $\mu$ is a   braiding map on  $\Vccgl$.
\end{lem}

\begin{prf}
Recall Definition \ref{def27}.
First,   note that, due to Remark \ref{rem_useful}, $\mu $ possesses  unitarity property \eqref{gen_br_uni}. Next, it is clear from \eqref{mu_form} that $\mu$ satisfies the singularity condition \eqref{gen_br_sing}. Hence, it remains to show that the quantum   Yang--Baxter equation \eqref{gen_br_ybe}   holds. Its proof employs the alternative expression for the action of $\mu$,
\begin{align}
 &\mu(z,x)\Big(T_{[m]}^{+24}(v)\ts \RR_{nm}^{  12}(u|v|2x+z ) \ts T_{[n]}^{+13}(u) \ts\R_{nm}^{  12}(u|v|z) (\vac\ot\vac)  \Big) \non \\
=\, &     \R_{nm}^{  12}(u|v|z)\ts  T_{[n]}^{+13}(u)  \ts \RR_{nm}^{  12}(u|v|2x+z )\ts  T_{[m]}^{+24}(v)  (\vac\ot\vac)
  ,\label{mybe1}
\end{align}
which is obtained from \eqref{mu_form} by the use of crossing symmetry  \eqref{RCSYMAx}. It is sufficient to check that the images of
\begin{align}
&T_{[k]}^{+36}(w)\ts \RR_{mk}^{23} \ts
T_{[m]}^{+25}(v)\ts \R_{mk}^{23} \ts
\RR_{nm}^{12}\ts \RR_{nk}^{13} \ts
T_{[n]}^{+14}(u)\ts \R_{nk}^{13} \ts
\R_{nm}^{12}(\vac\ot\vac\ot\vac)   \label{mybe2}
\end{align}
under   both sides of \eqref{gen_br_ybe} coincide.
In \eqref{mybe2},  we use the notation 
\begin{align*}
& \R_{nm}^{12}=\R_{nm}^{12}(u|v|z_1),\qquad
\RR_{nm}^{12}=\RR_{nm}^{12}(u|v|2x+2z_2+z_1),\qquad
 \R_{nk}^{13}= \R_{nk}^{13}(u|w| z_1+z_2),\\
& \RR_{nk}^{13}=\RR_{nk}^{13}(u|w|2x+z_1+z_2), \qquad
\R_{mk}^{23}=\R_{mk}^{23}(v|w|z_2),\qquad
\RR_{mk}^{23}= \RR_{mk}^{23}(v|w|2x+z_2) 
\end{align*}
and  the meaning of superscripts $1$--$6$ is analogous to \eqref{tfactors_sup} (except that the three copies of $\Vc$ are replaced by $\Vccgl$). 

Let us consider the action of the left-hand side of \eqref{gen_br_ybe} on \eqref{mybe2}. By applying $\mu_{23}(z_2,x)$ to \eqref{mybe2} and using the fact that
$T_{[k]}^{+36}(w)$  and
$\RR_{nm}^{12} $ commute,
  we get
\begin{align}
&\R_{mk}^{23} \ts T_{[m]}^{+25}(v)\ts 
\RR_{mk}^{23} \ts 
\RR_{nm}^{12} \ts T_{[k]}^{+36}(w)\ts 
\RR_{nk}^{13} \ts
T_{[n]}^{+14}(u)\ts \R_{nk}^{13} \ts
\R_{nm}^{12} (\vac\ot\vac\ot\vac).\label{mybe3}
\end{align}
Next, we apply $\mu_{13}(z_1+z_2,x)$ to \eqref{mybe3}, thus getting
 \begin{align}
&\R_{mk}^{23} \ts T_{[m]}^{+25}(v)\ts 
\RR_{mk}^{23} \ts 
\RR_{nm}^{12}\ts
\R_{nk}^{13} \ts
T_{[n]}^{+14}(u)\ts 
\RR_{nk}^{13} \ts
T_{[k]}^{+36}(w)\ts 
\R_{nm}^{12} (\vac\ot\vac\ot\vac).\label{mybe4}
\end{align}

By using the quantum Yang--Baxter equation \eqref{RYBE} one easily verifies the identities
\begin{align*}
\RR_{mk}^{23} \ts 
\RR_{nm}^{12} \ts
\R_{nk}^{13} 
=
\R_{nk}^{13} \ts
\RR_{nm}^{12} \ts
\RR_{mk}^{23}\fand
 \RR_{mk}^{23} \ts \RR_{nk}^{13} \ts \R_{nm}^{12}
=\R_{nm}^{12}
\ts \RR_{nk}^{13}\ts
\RR_{mk}^{23} .
\end{align*}
By employing these two equalities, along with commutation relations
\begin{gather*}
\RR_{mk}^{23}\ts T_{[n]}^{+14}(u)=T_{[n]}^{+14}(u)\ts \RR_{mk}^{23},\qquad
T_{[k]}^{+36}(w)\ts \R_{nm}^{12}=\R_{nm}^{12}\ts T_{[k]}^{+36}(w),\\
T_{[m]}^{+25}(v)\ts \R_{nk}^{13}= \R_{nk}^{13}\ts T_{[m]}^{+25}(v) ,
\end{gather*}
 we rewrite \eqref{mybe4} as
 \begin{align}
&\R_{mk}^{23} \ts 
\R_{nk}^{13} \ts
T_{[m]}^{+25}(v)\ts 
\RR_{nm}^{12}\ts
T_{[n]}^{+14}(u)\ts
 \R_{nm}^{12} \ts
\RR_{nk}^{13} \ts
\RR_{mk}^{23} \ts 
T_{[k]}^{+36}(w)(\vac\ot\vac\ot\vac)
.\label{mybe5}
\end{align}
Finally,  applying $\mu_{12}(z_1,x+z_2)$ to \eqref{mybe5} we obtain 
 \begin{align}
&\R_{mk}^{23} \ts 
\R_{nk}^{13} \ts
 \R_{nm}^{12} \ts
T_{[n]}^{+14}(u)\ts
\RR_{nm}^{12}\ts
T_{[m]}^{+25}(v)\ts 
\RR_{nk}^{13} \ts
\RR_{mk}^{23} \ts 
T_{[k]}^{+36}(w)(\vac\ot\vac\ot\vac)
.\label{mybe6}
\end{align}

Let us consider the action of the right-hand side of \eqref{gen_br_ybe} on \eqref{mybe2}. First, we observe that the quantum Yang--Baxter equation \eqref{RYBE} implies
$$
 \R_{mk}^{23} \ts
\RR_{nm}^{12}\ts \RR_{nk}^{13}=\RR_{nk}^{13}\ts
\RR_{nm}^{12}\ts  \R_{mk}^{23}\fand
\R_{mk}^{23} \ts\R_{nk}^{13} \ts
\R_{nm}^{12}=
\R_{nm}^{12}\ts\R_{nk}^{13} \ts \R_{mk}^{23} .
$$
By using these equalities and the commutation relations
$$
T_{[m]}^{+25}(v)\ts\RR_{nk}^{13}=\RR_{nk}^{13}\ts T_{[m]}^{+25}(v)
\fand
\R_{mk}^{23} \ts T_{[n]}^{+14}(u)=
T_{[n]}^{+14}(u) \ts \R_{mk}^{23} ,
$$
we rewrite \eqref{mybe2} as
\begin{align}
&T_{[k]}^{+36}(w)\ts 
\RR_{mk}^{23} \ts
\RR_{nk}^{13} \ts
T_{[m]}^{+25}(v)\ts 
\RR_{nm}^{12}\ts 
T_{[n]}^{+14}(u)\ts 
\R_{nm}^{12} \ts
\R_{nk}^{13} \ts
 \R_{mk}^{23}(\vac\ot\vac\ot\vac).\label{mybe7}
\end{align}
Next, we apply $\mu_{12}(z_1,x+z_2)$ to \eqref{mybe5}, thus getting 
\begin{align}
&T_{[k]}^{+36}(w)\ts 
\RR_{mk}^{23} \ts
\RR_{nk}^{13} \ts
\R_{nm}^{12} \ts
T_{[n]}^{+14}(u)\ts 
\RR_{nm}^{12}\ts 
T_{[m]}^{+25}(v)\ts 
\R_{nk}^{13} \ts
 \R_{mk}^{23}(\vac\ot\vac\ot\vac).\label{mybe9}
\end{align}

We shall also need the following consequences of  the   Yang--Baxter equation \eqref{RYBE}:
$$
\RR_{mk}^{23} \ts
\RR_{nk}^{13} \ts
\R_{nm}^{12}=\R_{nm}^{12}\ts
\RR_{nk}^{13} \ts \RR_{mk}^{23}\fand
 \RR_{mk}^{23}\ts \RR_{nm}^{12}\ts\R_{nk}^{13} 
=
\R_{nk}^{13} \ts \RR_{nm}^{12}\ts \RR_{mk}^{23}.
$$
By combining them with relations
\begin{gather*}
 \RR_{mk}^{23}\ts T_{[n]}^{+14}(u)=T_{[n]}^{+14}(u)\ts \RR_{mk}^{23},\qquad
T_{[m]}^{+25}(v)\ts 
\R_{nk}^{13}=\R_{nk}^{13}\ts T_{[m]}^{+25}(v), \\
T_{[k]}^{+36}(w)\ts 
\R_{nm}^{12}=\R_{nm}^{12}\ts T_{[k]}^{+36}(w) ,
\end{gather*}
we turn the expression in \eqref{mybe9} to
\begin{align}
\R_{nm}^{12} \ts
T_{[k]}^{+36}(w)\ts 
\RR_{nk}^{13} \ts
T_{[n]}^{+14}(u)\ts 
\R_{nk}^{13} \ts
\RR_{nm}^{12}\ts 
\RR_{mk}^{23} \ts
T_{[m]}^{+25}(v)\ts 
 \R_{mk}^{23}(\vac\ot\vac\ot\vac).\label{mybea}
\end{align}
By applying $\mu_{13}(z_1+z_2,x)$ to \eqref{mybea}, we get
\begin{align}
\R_{nm}^{12} \ts
\R_{nk}^{13} \ts
T_{[n]}^{+14}(u)\ts 
\RR_{nk}^{13} \ts
T_{[k]}^{+36}(w)\ts 
\RR_{nm}^{12}\ts 
\RR_{mk}^{23} \ts
T_{[m]}^{+25}(v)\ts 
 \R_{mk}^{23}(\vac\ot\vac\ot\vac).\label{mybeb}
\end{align}
As $T_{[k]}^{+36}(w)$ and 
$\RR_{nm}^{12}$ commute, the application of $\mu_{23}(z_2,x)$ to \eqref{mybeb} gives us
\begin{align}
\R_{nm}^{12} \ts
\R_{nk}^{13} \ts
T_{[n]}^{+14}(u)\ts 
\RR_{nk}^{13} \ts
\RR_{nm}^{12}\ts 
 \R_{mk}^{23}\ts 
T_{[m]}^{+25}(v)\ts
\RR_{mk}^{23} \ts
T_{[k]}^{+36}(w)(\vac\ot\vac\ot\vac)
.\label{mybec}
\end{align}

 The quantum Yang--Baxter equation \eqref{RYBE} implies
$$\RR_{nk}^{13} \ts
\RR_{nm}^{12}\ts 
 \R_{mk}^{23}=
\R_{mk}^{23} \ts
\RR_{nm}^{12}\ts \RR_{nk}^{13}
\fand
\R_{nm}^{12}\ts\R_{nk}^{13} \ts \R_{mk}^{23}
=\R_{mk}^{23} \ts\R_{nk}^{13} \ts
\R_{nm}^{12}.
$$
Finally, by employing these two equalities and the commutativity relations
$$
T_{[n]}^{+14}(u)\ts 
 \R_{mk}^{23}=
 \R_{mk}^{23}\ts T_{[n]}^{+14}(u)
\fand
\RR_{nk}^{13} \ts
T_{[m]}^{+25}(v)
=
T_{[m]}^{+25}(v)\ts\RR_{nk}^{13}, 
$$
we can rewrite \eqref{mybec} so that we get the expression in \eqref{mybe6}, which equals   the action of the left-hand side of \eqref{gen_br_ybe} on \eqref{mybe2}, thus concluding the proof.
\end{prf}

Following the second identity in \eqref{main_gamma},
we define the map $\nu$ by
$$
\nu(z,x)=\rho(z)\ts \gamma(z,x).
$$
Lemmas \ref{lem_rho} and \ref{lem_gamma} imply that its image belongs to
$ \Vccgl^{\ot 2}\ot\CC[z^{\pm 1}][[x^{-1},h]]$ and one can derive from these lemmas the explicit expression for the action of $\nu$,
 \begin{align}
 \nu(z,x)\Big( T_{[n]}^{+13}(u) \ts T_{[m]}^{+24}(v)(\vac\ot\vac ) \Big)  
= \,&\RR_{nm}^{  12}(u|v|2x+z-hN)^{-1}\cdotrl\left(   T_{[n]}^{+13}(u) \right.\label{nu_form}\\ 
  & \Big.\times\R_{nm}^{  12}(u|v| z-hN )^{-1}\ts  T_{[m]}^{+24}(v) \ts
\R_{nm}^{  12}(u|v|z)   \Big)(\vac\ot\vac ),\non
\end{align}
where ``$\cdotrl$'' denotes
the standard multiplication in 
$ ((\ndo\CC^N)^{\text{op}})^{\ot n} \ot (\ndo\CC^N)^{\ot m}$.

\begin{pro}\label{gen_thm_compatiblity}
The maps $ \mu$ and $\nu$ form a compatible pair on $\Vccgl$. Moreover, the pairs $(\mu,\nu)$ and $(\sigma,\rho)$ of $\Vccgl$ are equivalent.
\end{pro}

\begin{prf}
First, we observe that $\nu$ is an   intertwining map, so that $\mu$ is compatible; recall Definition \ref{def29}. Indeed,   property \eqref{mu1} holds due to Remark \ref{rem_useful}, while \eqref{mu2} is evident from   the expressions   \eqref{Y_form} and \eqref{nu_form} for the action of   $Y $ and $\nu$. In fact, using these formulae one easily derives the explicit formula for the action of   
$Y^\nu $,
 \begin{align}
 &Y^\nu(z,x)\Big( T_{[n]}^{+13}(u) \ts T_{[m]}^{+24}(v)(\vac\ot\vac ) \Big)\non\\  
&\qquad =  \RR_{nm}^{  12}(u|v|2x+z-hN)^{-1}\cdotrl\left(   T_{[n]}^{+13}(u|z)\ts    T_{[m]}^{+23}(v)   \right) \vac ,\label{Ynu-formula}
\end{align}
where ``$\cdotrl$'' has the same meaning as in \eqref{nu_form}. Hence, $(\mu,\nu)$ is a compatible pair.
As the  shift condition \eqref{shift_gamma}  is established in Lemma \ref{lem_gamma} and both identities in \eqref{main_gamma} clearly hold,   we conclude that the pairs $(\mu,\nu)$ and $(\sigma,\rho)$ are equivalent.
\end{prf}

Finally, we observe that the map $\gamma$ satisfies the following {\em hexagon relation}.

\begin{pro}
We have
\beq\label{hex_gamma} 
\gamma(z_1,x)\left(Y^\rho(z_2)\ot 1\right)=\left(Y^\rho(z_2)\ot 1\right)  \gamma_{13}(z_1+z_2,x)\ts \gamma_{23}(z_1,x).
\eeq
\end{pro}

\begin{prf}
It is sufficient to check that the images of \eqref{the_expression} 
under both sides of \eqref{hex_gamma} 
coincide. Of course, the superscripts $4$--$6$ in \eqref{the_expression} now represent the corresponding copies of $\Vccgl$; cf. \eqref{tfactors_sup}.
Regarding the left-hand side, by using   \eqref{Yrhoo} and \eqref{gamma_form}, we find that   the image of   \eqref{the_expression}  under
$\gamma(z_1,x)\left(Y^\rho(z_2)\ot 1\right)$ equals
\beq\label{image_1}
T_{[n]}^{+14 }(u|z_2) \ts T_{[m]}^{+24 }(v)\ts 
\RR_{mk}^{23}(v|w|2x+z_1-hN)^{-1}\ts
\RR_{nk}^{13}(u|w|2x+z_1+z_2-hN)^{-1}\ts
T_{[k]}^{+35}(w)(\vac\ot\vac ).
\eeq
On the other hand, by \eqref{gamma_form},
the image of  \eqref{the_expression}  under
$\gamma_{13}(z_1+z_2,x)  \gamma_{23}(z_1,x)$ is equal  to
$$
T_{[n]}^{+14 }(u ) \ts T_{[m]}^{+25 }(v)\ts 
\RR_{mk}^{23}(v|w|2x+z_1-hN)^{-1}\ts
\RR_{nk}^{13}(u|w|2x+z_1+z_2-hN)^{-1}\ts
T_{[k]}^{+36}(w)(\vac\ot\vac\ot\vac ).
$$
Finally,   \eqref{Yrhoo} implies that $Y^\rho(\cdot, z_2)\ot 1$
maps this expression  to \eqref{image_1}, as required.
\end{prf}

\subsection{Action of the  braiding   maps  on the quantum Feigin--Frenkel center}\label{section_04}
In this subsection, we investigate   fixed points of the compatible pairs $(\sigma,\rho)$ and $(\mu,\nu)$. Our key results are Lemmas \ref{fixed_lemma} and \ref{lemma_fixed_mu}, which we utilize in Sections \ref{section_06} and \ref{section_07} below.
From now on, we shall write bar over the number  in the superscript to indicate that the corresponding $R$-matrix product  appears in the reversed order, e.g., 

\begin{align*}
&\R_{nm}^{\bar{1}2}(u|v|z)= \prod_{i=1,\dots,n}^{\longleftarrow}
\prod_{j=n+1,\ldots,n+m}^{\longleftarrow} \R_{ij} ,\qquad
 &\R_{nm}^{1\bar{2}}(u|v|z)= \prod_{i=1,\dots,n}^{\longrightarrow}
\prod_{j=n+1,\ldots,n+m}^{\longrightarrow} \R_{ij},
\\	
&\R_{nm}^{\bar{1}\bar{2}}(u|v|z)= \prod_{i=1,\dots,n}^{\longleftarrow}
\prod_{j=n+1,\ldots,n+m}^{\longrightarrow} \R_{ij},&\text{where }
\R_{ij}=\R_{ij}(z+u_i -v_{j-n});
\end{align*}
recall \eqref{rnm12}. 
We apply the same convention to the $R$-matrix products defined by \eqref{rnm123}.
Also, we denote the reversed   products of the operator matrices from \eqref{teen} by
$$
\cev{T}^{\pm 12}_{[n]}(u)=T_n^\pm(u_n )\ldots T_1^\pm(u_1)\fand 
\cev{T}^{\pm 12}_{[n]} (u|z)=T^\pm_n(z+u_n)\ldots T^\pm_1(z+u_1).
$$
	Finally, for any $m$ and $n$   we shall consider the tensor products of the series \eqref{tm},
	$$ \tmn= \tm(u)\ot \tn(v)\in\Vcc\ot\Vcc[[u,v]].$$ 

	Consider the maps $\sigma $ and $\rho$  introduced in   Lemmas \ref{lem_sigma} and \ref{lem_rho}, respectively.

\begin{lem}\label{fixed_lemma}
\upshape{(a)}   The following identities hold in $\Vcc$:
\beq\label{fixed_rho}
\sigma(z)\Big(\tmn\Big)
 =\rho(z )\Big(\tmn\Big)=\tmn\quad\text{for any }m\text{ and }n.
\eeq
\upshape{(b)} The following identities hold in $\Vccgl$:
\beq\label{fixed_rho_2}
\sigma(z )\Big(w\ot\tN(u) \Big)=
\rho(z )\Big(w\ot\tN(u)\Big)=
w\ot\tN(u)\quad\text{for any }w\in\Vccgl.
\eeq
\end{lem}

\begin{prf}
Before we start with the proof, we introduce some new notation. Let
$$
\R_{nm}^{12}=\R_{nm}^{12}(u_{[n]}|v_{[m]}|z)\fand \mathcal{E}_{[n,m]}=\mathcal{E}_{[n ]}\ot \mathcal{E}_{[ m]},
$$
where $u_{[n]}$ (resp. $\mathcal{E}_{[n ]}$) are given by \eqref{uem0} and \eqref{uem1} (resp. \eqref{fusionA} and \eqref{fusionBCD}), depending on the type of $\g_N$. In the proof, we shall use the following well-known properties of the idempotent map $\mathcal{E}_{[n ]}$ (see, e.g.,  \cite[Lemma 2.3]{BJK} and the proof of \cite[Thm. 2.4]{JKMY}):
\beq\label{ysm-map}
\mathcal{E}^1_{[n ]}\ts \R_{nm}^{12}=\R_{nm}^{\bar{1}2}\ts \mathcal{E}^1_{[n ]},\quad
\mathcal{E}^2_{[m ]}\ts \R_{nm}^{12}=\R_{nm}^{1\bar{2}}\ts \mathcal{E}^2_{[m ]},\quad
\mathcal{E}^1_{[n ]}\ts T^{+ 12}_{[n]} (u_{[n]})=\cev{T}^{+ 12}_{[n]} (u_{[n]})\ts \mathcal{E}^1_{[n ]}.
\eeq
Note that the above  equalities  imply
\begin{align}
\mathcal{E}^{12}_{[n,m]}\ts \R_{nm}^{12}&=\R_{nm}^{\bar{1}\bar{2}}\ts \mathcal{E}^{12}_{[n,m]},\label{ysm-map-2}\\ 
\mathcal{E}^{12}_{[n,m]}\ts T^{+ 13}_{[n]} (u_{[n]})\ts T^{+ 24}_{[m]} (v_{[m]})&=
 \cev{T}^{+ 24}_{[m]} (v_{[m]})\ts \cev{T}^{+ 13}_{[n]} (u_{[n]})\ts  \mathcal{E}^{12}_{[n,m]}.\label{ysm-map-2-2}
\end{align}

\noindent\upshape{(a)}   Let us prove that all coefficients of $\tmn$ are fixed points of $\sigma$. By \eqref{sigma_formula}, the image of $\tmn$ under the braiding map $\sigma(z)$ equals
$$
\tr\ts 
\mathcal{E}^{12}_{[n,m]}\ts
\R_{nm}^{12}\ts
T^{+ 13}_{[n]} (u_{[n]})\ts 
T^{+ 24}_{[m]} (v_{[m]})\ts
(\R_{nm}^{12})^{-1}(\vac\ot\vac ),
$$
where the trace is taken over the tensor factors $1,\ldots, n+m$. By employing the identities     \eqref{ysm-map-2} and \eqref{ysm-map-2-2}, along with $(\mathcal{E}^{12}_{[n,m]})^2=\mathcal{E}^{12}_{[n,m]}$, we rewrite the given expression 
as follows:
\begin{align*}
&\tr\ts
\R_{nm}^{\bar{1}\bar{2}}\ts
\mathcal{E}^{12}_{[n,m]}\ts
T^{+ 13}_{[n]} (u_{[n]})\ts 
T^{+ 24}_{[m]} (v_{[m]})\ts
(\R_{nm}^{12})^{-1}(\vac\ot\vac )\\
=&
\tr\ts
\R_{nm}^{\bar{1}\bar{2}}\ts
\mathcal{E}^{12}_{[n,m]}\ts
\mathcal{E}^{12}_{[n,m]}\ts
T^{+ 13}_{[n]} (u_{[n]})\ts 
T^{+ 24}_{[m]} (v_{[m]})\ts
(\R_{nm}^{12})^{-1}(\vac\ot\vac )\\
=&
\tr\ts
\mathcal{E}^{12}_{[n,m]}\ts
\R_{nm}^{12}\ts
\cev{T}^{+ 13}_{[n]} (u_{[n]})\ts 
\cev{T}^{+ 24}_{[m]} (v_{[m]})\ts
(\R_{nm}^{\bar{1}\bar{2}})^{-1}
\mathcal{E}^{12}_{[n,m]}(\vac\ot\vac ).
\end{align*}
Next, we use the cyclic property of the trace to move the left copy of $\mathcal{E}^{12}_{[n,m]}$ to the right:
\begin{align*}
&\tr\ts
\R_{nm}^{12}\ts
\cev{T}^{+ 13}_{[n]} (u_{[n]})\ts 
\cev{T}^{+ 24}_{[m]} (v_{[m]})\ts
(\R_{nm}^{\bar{1}\bar{2}})^{-1}
(\mathcal{E}^{12}_{[n,m]})^2(\vac\ot\vac )\\
 =&
\tr\R_{nm}^{12}\ts
\cev{T}^{+ 13}_{[n]} (u_{[n]})\ts 
\cev{T}^{+ 24}_{[m]} (v_{[m]})\ts
(\R_{nm}^{\bar{1}\bar{2}})^{-1}
\mathcal{E}^{12}_{[n,m]}(\vac\ot\vac ).
\end{align*}
Finally, we employ the equalities   \eqref{ysm-map-2} and \eqref{ysm-map-2-2}  once again, to move $\mathcal{E}^{12}_{[n,m]}$ to the left:
$$
\tr\R_{nm}^{12}\ts
\mathcal{E}^{12}_{[n,m]}\ts
T^{+ 13}_{[n]} (u_{[n]})\ts 
T^{+ 24}_{[m]} (v_{[m]})\ts
(\R_{nm}^{12})^{-1}(\vac\ot\vac )
.
$$
By the cyclic property of the trace, the terms $\R_{nm}^{12}$ and $(\R_{nm}^{12})^{-1}$ cancel, so   we get
$$
\tr 
\mathcal{E}^{12}_{[n,m]}\ts
T^{+ 13}_{[n]} (u_{[n]})\ts 
T^{+ 24}_{[m]} (v_{[m]})(\vac\ot\vac )=\tmn,
$$
as required.

The proof that all coefficients of $\tmn$ are fixed points of $\rho$ is similar. However,   some of its intermediate steps differ, so we present it in detail.
Let 
$$
\R_{nm}^{12\ts *}=\R_{nm}^{12}(u_{[n]}|v_{[m]}|z+h\varepsilon c_{\text{crit}})^{-1}
=\begin{cases}
\R_{nm}^{12}(u_{[n]}|v_{[m]}|z-hN)^{-1}&\text{ for }\g_N=\gl_N,\\
\R_{nm}^{12}(u_{[n]}|v_{[m]}|z-2h\kp)^{-1}&\text{ for }\g_N=\on_N,\spn_N.
\end{cases}
$$
First, we remark that the first two formulae in \eqref{ysm-map} hold for $\R_{nm}^{12\ts *} $ as well.
By \eqref{rho_formula}, the image of $\tmn$ under the   map $\rho(z)$ equals
$$
\tr\ts 
\mathcal{E}^1_{[n ]}  \ts
\mathcal{E}^2_{[ m]}\ts
T^{+ 13}_{[n]} (u_{[n]})\ts 
 \R_{nm}^{12\ts*}  \ts
T^{+ 24}_{[m]} (v_{[m]})\ts
 \R_{nm}^{12}(\vac\ot\vac ) ,
$$
where the trace is again taken over the tensor factors $1,\ldots, n+m$.
We now use
the  identities in \eqref{ysm-map} along with $(\mathcal{E}^1_{[n ]})^2=\mathcal{E}^1_{[n ]}$ and $(\mathcal{E}^2_{[ m]})^2=\mathcal{E}^2_{[ m]}$ to rewrite the given expression as follows:
\begin{align*}
&\tr\ts 
\cev{T}^{+  13}_{[n]} (u_{[n]})\ts 
\mathcal{E}^1_{[n ]}  \ts
 \R_{nm}^{1\bar{2}\ts*}  \ts
\mathcal{E}^2_{[ m]}\ts
T^{+ 24}_{[m]} (v_{[m]})\ts
 \R_{nm}^{12} (\vac\ot\vac )\\
=&
\tr\ts
\cev{T}^{+ 13 }_{[n]} (u_{[n]})\ts 
\mathcal{E}^1_{[n ]}  \ts \mathcal{E}^1_{[n ]}  \ts
 \R_{nm}^{1\bar{2}\ts*}  \ts
\mathcal{E}^2_{[ m]}\ts \mathcal{E}^2_{[ m]}\ts
T^{+ 24}_{[m]} (v_{[m]})\ts
 \R_{nm}^{12} (\vac\ot\vac ).
\end{align*}
Next, we use the same identities to move the left (resp. right) copies of
$ \mathcal{E}^1_{[n ]}$ and $ \mathcal{E}^2_{[ m]}$ all the way to the left (resp. right), thus getting
\begin{align*}
&
\tr\ts
\mathcal{E}^1_{[n ]}  \ts
\mathcal{E}^2_{[ m]}\ts
T^{+ 13}_{[n]} (u_{[n]})\ts 
 \R_{nm}^{\bar{1}2\ts*}  \ts
\cev{T}^{+ 24}_{[m]} (v_{[m]})\ts
 \R_{nm}^{\bar{1}\bar{2}} \ts
 \mathcal{E}^1_{[n ]}  \ts
 \mathcal{E}^2_{[ m]}(\vac\ot\vac ).
\end{align*}
Using the cyclic property of the trace we move the left copies of $ \mathcal{E}^1_{[n ]}$ and $ \mathcal{E}^2_{[ m]}$ to the right:
\begin{align*}
&
\tr\ts
T^{+ 13}_{[n]} (u_{[n]})\ts 
 \R_{nm}^{\bar{1}2\ts*}  \ts
\cev{T}^{+ 24}_{[m]} (v_{[m]})\ts
 \R_{nm}^{\bar{1}\bar{2}} \ts
 \mathcal{E}^1_{[n ]}  \ts
 \mathcal{E}^2_{[ m]}(\vac\ot\vac ).
\end{align*}
Finally, we use the equalities in \eqref{ysm-map} once again, to move $ \mathcal{E}^1_{[n ]}$ and $ \mathcal{E}^2_{[ m]}$ to the left:
\begin{align*}
&
\tr\ts
T^{+ 13}_{[n]} (u_{[n]})\ts 
 \mathcal{E}^1_{[n ]}  \ts
 \R_{nm}^{12\ts*}  \ts
 \mathcal{E}^2_{[ m]} \ts
T^{+ 24}_{[m]} (v_{[m]})\ts
 \R_{nm}^{12} (\vac\ot\vac )
.
\end{align*}
By the cyclic property of the trace, this can be written equivalently as
\beq\label{RCSYMBCDxx-7}
\tr\ts
T^{+ 13}_{[n]} (u_{[n]})\ts 
 \mathcal{E}^1_{[n ]}  \ts
 \mathcal{E}^2_{[ m]} \ts
T^{+ 24}_{[m]} (v_{[m]}) 
\left(
 \R_{nm}^{12\ts*}  \cdotlr
 \R_{nm}^{12} \right)(\vac\ot\vac )
,
\eeq
where the symbol ``$\cdotlr$'' indicates the standard multiplication in the tensor product algebra
$ (\ndo\CC^N)^{\ot n} \ot ((\ndo\CC^N)^{\text{op}})^{\ot m}$.
By the crossing symmetry properties \eqref{RCSYMAx} and \eqref{RCSYMBCDxx}, we have
$ \R_{nm}^{12\ts*}  \cdotlr
 \R_{nm}^{12} =1$, so the expression in \eqref{RCSYMBCDxx-7} is equal to
\begin{align*}
&\tr\ts
T^{+ 13}_{[n]} (u_{[n]})\ts 
 \mathcal{E}^1_{[n ]}  \ts
 \mathcal{E}^2_{[ m]} \ts
T^{+ 24}_{[m]} (v_{[m]})(\vac\ot\vac ) \\
=&
\tr\ts
 \mathcal{E}^1_{[n ]}  \ts
T^{+ 13}_{[n]} (u_{[n]})\ts 
 \mathcal{E}^2_{[ m]} \ts
T^{+ 24}_{[m]} (v_{[m]})(\vac\ot\vac )
=\tmn,
\end{align*}
as required.

\noindent\upshape{(b)} Both equalities can be verified by arguing  as in the proof of the first assertion of   lemma. 
In addition, the analogous arguments can be found in the proof of Lemma \ref{lemma_fixed_mu} below, so we omit the proof.
\end{prf}

\begin{rem}\label{qdet}
The series $\tN(u)\in \Vccgl[[u]]$, which appears in \eqref{fixed_rho_2}, is   the well-known {\em quantum determinant} $\qdet T^+(u)$ for the dual Yangian $\Ydgl$.
Moreover, the corresponding property of the braiding map,
$\sigma(z ) (w\ot\tN(u)  )=w\ot\tN(u)$ 
actually holds in $\Vcgl$ for any $c\in\CC$, but we   use it later   at the critical level only.
\end{rem}

The first statement of the next corollary for $\mathfrak{g}_N=\mathfrak{o}_N,\mathfrak{sp}_N$ was   established in \cite[Thm. 2.5]{BJK}.
Although we obtain it here as an immediate consequence of \eqref{mu1} and Lemma \ref{fixed_lemma}, we believe   it is worth to single it out, as the braiding of an arbitrary $h$-adic quantum vertex algebra does not need to exhibit such a property; see \cite[Thm. 1.4]{DGK}. The second statement follows by combining Theorem \ref{thm_ff} and the first statement.

\begin{kor} 
 \upshape{(a)} For any $m$ and $n$ we have
$$
\Sc(z) \Big(\tmn\Big)=\tmn.
$$
 \upshape{(b)}
For $\mathfrak{g}_N$   of type $A,B$,
the restriction of    $\Sc  $ to $\mathfrak{z}(\Vcc)^{\ot 2}  $ 
is the identity.
\end{kor}

Analogously to the previous corollary, one also establishes its counterpart for $\sigma$.
\begin{kor}
For $\mathfrak{g}_N$  of type $A,B$, 
the restriction of   $\sigma $ to $\mathfrak{z}(\Vcc)^{\ot 2}  $ 
is the identity.
\end{kor}

Finally, we obtain a partial   analogue of Lemma \ref{fixed_lemma} for the   map $\mu$  defined by  \eqref{mu_form_second}.

\begin{lem}\label{lemma_fixed_mu} 
The  braiding map $\mu$ satisfies
\begin{align}
&\mu(z,x )\Big(\tjj\Big)=\tjj ,\label{fixed_mu}\\
&\mu(z,x )\Big(w\ot \tN(u) \Big)=w\ot \tN(u)\text{ for all }w\in\Vccgl.\label{fixed_mu_two}
\end{align}
 \end{lem}

\begin{prf}
Let us prove   \eqref{fixed_mu}. We shall use the notation
\begin{align*}
\R_{12}=\R_{12}(z+u-v),\,\,
\RR_{12}=\RR_{12}(2x+z+u+v),\,\,
\RR^*_{12}=\RR_{12}(2x+z+u+v-hN)^{-1}.
\end{align*}
By \eqref{mu_form}, the image of $\tjj$ under the map $\mu$ is equal to
$$
\tr\ts 
\RR^*_{12}\cdotlr\left(
\R_{12}\ts
T_{13}^+(u)\ts
\RR_{12}\ts
T_{24}^+(v)\ts
(\R_{12})^{-1}
\right)(\vac\ot\vac ),
$$
where the trace is taken over the tensor factors $1$ and $2$. Due to  the cyclic property of the trace and the commutation relation $(\R_{12})^{-1} \RR^*_{12}=\RR^*_{12}(\R_{12})^{-1}$, this equals
$$
\tr\ts 
\RR^*_{12}\cdotlr\left(
(\R_{12})^{-1}\ts
\R_{12}\ts
T_{13}^+(u)\ts
\RR_{12}\ts
T_{24}^+(v)\ts
\right)(\vac\ot\vac )
=
\tr\ts 
\RR^*_{12}\cdotlr\left(
T_{13}^+(u)\ts
\RR_{12}\ts
T_{24}^+(v)\ts
\right)(\vac\ot\vac )
.
$$
We use   the cyclic property of the trace once more to rewrite the given expression as
$$
\tr\ts 
\RR^*_{12}\cdotrl\left(
T_{13}^+(u)\ts
\RR_{12}\ts
T_{24}^+(v)
\right)(\vac\ot\vac )
=
\tr\ts 
T_{13}^+(u)
\left(\RR^*_{12}\cdotrl\RR_{12}\right)
T_{24}^+(v)(\vac\ot\vac )
.
$$
It remains to observe that the expression above is equal to $\tjj$ since the  unitarity and the crossing symmetry relations \eqref{Runi} and \eqref{RCSYMAx} imply 
$\RR^*_{12}\cdotrl\RR_{12}=1$.

To prove   \eqref{fixed_mu_two}, we   use  the  property of the anti-symmetrizer,
\beq\label{antisym}
\mathcal{E}_{[N]}\ts\R_{0N}(u+(N-1)h)\ldots \R_{02}(u+ h)\ts\R_{01}(u)=\mathcal{E}_{[N]},
\eeq
which holds in the case $\mathfrak{g}_N=\mathfrak{gl}_N$; see \cite[Ch. 1]{M}.
Furthermore, we use the  notation
\begin{align*}
\R_{nN}^{12}=\R_{nN}^{12}(v|u_{[N]}|z),\,\,
\RR_{nN}^{12}=\RR_{nN}^{12}(v|u_{[N]}|2x+z),\,\,
\RR_{nN}^{12\ts*}=\RR_{nN}^{12}(v|u_{[N]}|2x+z-hN)^{-1},
\end{align*}
where the $N$-tuple $u_{[N]}$ is given by \eqref{uem0}.
By employing  \eqref{antisym}, one easily proves
\beq\label{antisymm}
\mathcal{E}_{[N]}\ts (\R_{nN}^{12})^{\pm 1}=
\mathcal{E}_{[N]}\ts \RR_{nN}^{12}=
\mathcal{E}_{[N]}\ts\RR_{nN}^{12\ts*}= \mathcal{E}_{[N]}.
\eeq
Due to \eqref{mu_form}, the image of $T_{[n]}^{+13}(v)\ot \tN(u)$ under the map  $\mu$ is equal to
\beq\label{mu_lemma_tmp-1}
\tr\ts
\mathcal{E}_{[N]}^2\ts
\RR_{nN}^{12\ts*}\cdotlr
\left(
\R_{nN}^{12}\ts 
T_{[n]}^{+13}(v)\ts
\RR_{nN}^{12}\ts
T_{[N]}^{+24}(u_{[N]})\ts
(\R_{nN}^{12})^{-1}
\right)(\vac\ot\vac ),
\eeq
where the trace is taken over the tensor factors $n+1,\ldots ,n+N$
and ``$\cdotlr$'' denotes the standard multiplication in $(\ndo\CC^N)^{\ot n}\ot ((\ndo\CC^N)^{\text{op}})^{\ot N}$.
By employing \eqref{antisymm} and the third equality in
\eqref{ysm-map} with $n=N$, we find that \eqref{mu_lemma_tmp-1} equals
\begin{align*}
&
\tr\ts
\mathcal{E}_{[N]}^2\ts
\RR_{nN}^{12\ts*}\cdotlr
\left(
\R_{nN}^{12}\ts 
T_{[n]}^{+13}(v)\ts
\RR_{nN}^{12}\ts
T_{[N]}^{+24}(u_{[N]})\ts
(\R_{nN}^{12})^{-1}
\right)(\vac\ot\vac )\\
=&
\tr\ts
\RR_{nN}^{12\ts*}\cdotlr
\left(
\mathcal{E}_{[N]}^2\ts
\R_{nN}^{12}\ts 
T_{[n]}^{+13}(v)\ts
\RR_{nN}^{12}\ts
T_{[N]}^{+24}(u_{[N]})\ts
(\R_{nN}^{12})^{-1}
\right)(\vac\ot\vac )\\
=&
\tr\ts
\RR_{nN}^{12\ts*}\cdotlr
\left(
T_{[n]}^{+13}(v)\ts
\mathcal{E}_{[N]}^2\ts 
\RR_{nN}^{12}\ts
T_{[N]}^{+24}(u_{[N]})\ts
(\R_{nN}^{12})^{-1}
\right)(\vac\ot\vac )\\
=&
\tr\ts
\RR_{nN}^{12\ts*}\cdotlr
\left(
T_{[n]}^{+13}(v)\ts
\cev{T}_{[N]}^{+24}(u_{[N]})\ts
\mathcal{E}_{[N]}^2\ts 
(\R_{nN}^{12})^{-1}
\right)(\vac\ot\vac )\\
=&
\tr\ts
\RR_{nN}^{12\ts*}\cdotlr
\left(
T_{[n]}^{+13}(v)\ts
\cev{T}_{[N]}^{+24}(u_{[N]})\ts
\mathcal{E}_{[N]}^2 
\right)(\vac\ot\vac ).
\end{align*}
By the cyclic property of the trace and the last equality in \eqref{antisymm}, this equals
$$
\tr\ts
\left(\mathcal{E}_{[N]}^2\ts
\RR_{nN}^{12\ts*}\right)\cdotlr
\left(
T_{[n]}^{+13}(v)\ts
\cev{T}_{[N]}^{+24}(u_{[N]})\ts 
\right)(\vac\ot\vac )
=
\tr\ts
T_{[n]}^{+13}(v)\ts
\cev{T}_{[N]}^{+24}(u_{[N]})\ts
 \mathcal{E}_{[N]}^2 (\vac\ot\vac ).
$$
Finally, we use the third equality in
\eqref{ysm-map} to move $\mathcal{E}_{[N]}$ to the left, thus getting
$$
\tr\ts
T_{[n]}^{+13}(v)\ts
 \mathcal{E}_{[N]}^2\ts 
T_{[N]}^{+24}(u_{[N]}) (\vac\ot\vac )=T_{[n]}^{+13}(v)\ot \tN(u),
$$
as required.
\end{prf}

\section{Generalized Yangians  as  \texorpdfstring{$(\sigma, \rho )$-deformed}{(sigma,rho)-deformed} \texorpdfstring{$\Vc$}{Vc(gN)}-modules}\label{section_06}

In this section, we associate $(\sigma,\rho)$-deformed modules with the generalized Yangians of types $A$--$D$ and demonstrate their applications.
Recall the $R$-matrix $R(u)$ defined by \eqref{RA} (resp. \eqref{RBCD}) for $\g_N=\gl_N$ (resp. $\g_N=\on_N,\spn_N$). We shall also need the $R$-matrix $\RRR(u)$, which is a polynomial in the variable $u$,  
\beq\label{Rhat-R-matrices}
\RRR(u)=\begin{cases}
u\ts R(u)&\text{ for }\g_N=\gl_N,\\
u(u-h\kp) R(u)&\text{ for }\g_N=\on_N,\spn_N.
\end{cases}
\eeq

Following  Krylov and Rybnikov\cite[Sect. 3.1]{KR}, we define the  algebra $ \Yg$      as the $h$-adically completed associative algebra over the ring $\CC[[h]]$ generated by the     elements $\lambda_{ij}^{(  r)}$, where $i,j=1,\ldots ,N$ and $r\in\ZZ,$ subject to the defining relations
\begin{align}
\RRR(u-v)\ts \Lc_1 (u)\ts\Lc_2(v) =\Lc_2(v)\ts\Lc_1 (u)\ts\RRR(u-v)\label{YY1}.
\end{align}
The generator matrix $\Lc(u)$    is defined by
\beq\label{YY2}
\Lc(u)=\sum_{i,j=1}^N e_{ij}\ot \lambda_{ij}(u),\quad\text{where}\quad \lambda_{ij}(u)=\delta_{ij}-h\sum_{r\in\mathbb{Z}} \lambda_{ij}^{(-r)}u^{r}.
\eeq
Consider the   two-sided ideals $I_p$ in $ \Yg$, where $p\geqslant 1$, generated by the elements
$$
\lambda_{ij}^{(-r)}\quad\text{with}\quad r\geqslant p\Fand h^{p }\cdot 1.
$$
Clearly, the family $(I_p)_{p\geqslant 0}$, where we set $I_0=\Yg$, satisfies the constraints in \eqref{top_submodules}. Denote by $\tau$ the topology on
$ \Yg$ defined by the basis $a+I_p$ for all $a\in \Ygl$ and $p\geqslant 0$. Consider the  algebra completed with respect to $\tau$,
$
 \lim_{\longleftarrow} \Yg/I_p.
$
For $\g_N=\gl_N$ (resp. $\g_N=\on_N,\spn_N$), we denote this completed algebra by $ \YYgl$ (resp. $\YYg^\prime$). Finally, for
$\g_N=\on_N,\spn_N$, we   define the algebra $\YYg$ as the quotient of 
 $\YYg^\prime$ over the $\tau$-completed two-sided ideal generated by the relations
\beq\label{YY3}
\Lc(u)\Lc(u+h\kp)^\prime=1.
\eeq
In contrast with  the original definition \cite[Sect. 3.1]{KR} for $\g_N=\gl_N$, which is given over the complex field, we consider algebras over the ring $\CC[[h]]$ so that they fit the $h$-adic  quantum vertex algebra setting.
We refer to the algebras $\YYg$ as {\em generalized Yangians}.

In analogy with \eqref{teen}, we introduce the notation
\beq\label{leen}
\Lc_{[n]}(u)=\Lc_1(u_1)\ldots \Lc_n(u_n )\fand \Lc_{[n]} (u|z)=\Lc_1(z+u_1)\ldots \Lc_n(z+u_n).
\eeq

\begin{thm}\label{generalized_yangian_thm}
For any $c\in\CC$, there exists a unique  structure of   $(\sigma,\rho )$-deformed $\Vc$-module  over $\YYg$ such that the module map is defined by
\beq\label{Y-module-map}
Y_{\YYg}\big(T_{[m]}^+ (u)\vac,z\big)=\Lc_{[m]}  (u|z). 
\eeq
\end{thm}

\begin{prf}
The map $Y_{\YYg}(\cdot ,z)$ is well-defined by \eqref{Y-module-map} due to the fact that the defining relations for the dual Yangian, \eqref{DY1} and \eqref{DY3}, both for the matrix $T^+(u)$, are of the same form as the defining relations \eqref{YY1} and \eqref{YY3} for the algebra $\YYg$. In particular, note that by replacing the $R$-matrix $R(u-v)$ in \eqref{DY1} by $\RRR(u-v)$, one obtains  equivalent relations.
Thus, to prove the theorem, it is sufficient to check that the requirements imposed by Definition \ref{defn_deformed_top} hold. It is clear from the definition of the topology $\tau$ that the image of $Y_{\YYg} (\cdot,z)$ belongs to $\YYg((z^{-1}))_{\tau}$ and, in addition, that this map satisfies the truncation condition \eqref{trunc_cond}. Moreover, the vacuum property \eqref{vac_cond_deformed}  obviously holds.

As for the $\rho$-associativity,   i.e. \eqref{rho_assoc_top} with $\nu=\rho$, one easily derives from \eqref{Yrhoo} and \eqref{Y-module-map} that both   sides of \eqref{rho_assoc_top} map the expression
$T_{[n]}^{+13}(u)  T_{[m]}^{+24}(v)$, with $n$ and $m$ being arbitrary and $u=(u_1,\ldots ,u_n)$, $v=(v_1,\ldots ,v_m)$, to
$$
\Lc_{1k}(z_2+z_0+u_1)\ldots
\Lc_{nk}(z_2+z_0+u_n)\ts
\Lc_{n+1\ts k}(z_2+v_1)\ldots
\Lc_{n+m\ts k}(z_2+v_m) 
$$
for $k=n+m+1$, so it is evident that   $Y_{\YYg}(\cdot,z)$ possesses the $\rho$-associativity property.

It remains to show that  $Y_{\YYg}(\cdot,z)$ satisfies the $\sigma$-commutativity, i.e. \eqref{sigma_loc_top}  with $\mu=\sigma$. By \eqref{sigma_formula} and \eqref{Y-module-map}, the left-hand side of \eqref{sigma_loc_top} maps the expression $T_{[n]}^{+13}(u)  T_{[m]}^{+24}(v)$ to
$$
\R_{nm}^{12}(u|v|-z_2+z_1)\ts \Lc_{[n]}^{13}(u|z_1)\ts
\Lc_{[m]}^{23}(v|z_2)\ts \R_{nm}^{12}(u|v|-z_2+z_1)^{-1}.
$$
By canceling the normalizing functions, this can be written via the $R$-matrices \eqref{Rhat-R-matrices} as
\beq\label{sigma-locality-lhs}
\RRR_{nm}^{12}(u|v|  -z_2+z_1)\ts \Lc_{[n]}^{13}(u|z_1)\ts
\Lc_{[m]}^{23}(v|z_2)\ts \RRR_{nm}^{12}(u|v| -z_2+z_1 )^{-1},
\eeq
where $\RRR_{nm}^{12}(u|v|z)$ is defined as in \eqref{rnm12} with the factors $\R_{ij}(z+u_i -v_{j-n})$ being replaced by $\RRR_{ij}(z+u_i -v_{j-n})$. On the other hand, the image of $T_{[n]}^{+13}(u)  T_{[m]}^{+24}(v)$ under the right-hand side of \eqref{sigma_loc_top} equals
\beq\label{sigma-locality-lhs-2}
\Lc_{[m]}^{23}(v|z_2)\ts \Lc_{[n]}^{13}(u|z_1).
\eeq
Finally, it is clear from  the defining relations \eqref{YY1} that \eqref{sigma-locality-lhs} and \eqref{sigma-locality-lhs-2} coincide, so we conclude that the map $Y_{\YYg}(\cdot,z)$ is $\sigma$-commutative.
\end{prf}
 
Next, we consider
 the structure of $(\sigma, \rho )$-deformed $\Vcc$-module over $\YYg$, i.e. we restrict the setting to the critical level.
The next corollary is a  direct consequence of the $\sigma$-commutativity \eqref{sigma_loc_top}  and the properties of the map $\sigma$ given by \eqref{fixed_rho} and \eqref{fixed_rho_2}.

\begin{kor}\label{koryy} 
\upshape{(a)} The coefficients of all
$$
\mathbb{L}_{[n]}  (z)\coloneqq Y_{\YYg} \big(\tm (0), z \big) \in \YYg[[z^{\pm 1}]]
$$
generate a commutative subalgebra in  $\YYg$.

\noindent\upshape{(b)} The coefficients of  
$$
\lN (z) = Y_{\YYgl} \big(\tN (0), z \big)\in \YYgl[[z^{\pm}]]
$$
belong to the center of   $\YYgl$.
\end{kor}

\begin{rem}
In the case $\mathfrak{g}_N=\mathfrak{gl}_N$, the fact that the coefficients of 
$\mathbb{L}_{[m]}  (z)$ generate a commutative subalgebra was originally established in \cite[Sect. 3]{KR}.  Moreover, a different vertex-operator theoretic  interpretation of the series $\mathbb{L}_{[m]}  (z)$ for $\mathfrak{g}_N=\mathfrak{gl}_N$  and a direct proof of the second assertion of Corollary \ref{koryy}  was given in \cite[Sect. 6, 7]{BK}.
\end{rem}

The constructions from this section can be    specialized to fit the     (ordinary) $h$-adically completed Yangian $\Y$. Following \eqref{tDY}, one    redefines the series $\lambda_{ij}(u)$ in \eqref{YY2} as
$
\lambda_{ij}(u)=\delta_{ij}+h\sum_{r\geqslant 1} \lambda_{ij}^{(r)}u^{-r}
$
and sets $\tau$ to be the $h$-adic topology. Then, the next   analogues of Theorem \ref{generalized_yangian_thm} and Corollary \ref{koryy}   follow by  simply   adjusting the above  arguments.

\begin{thm}\label{kor_tm_tm_Yang} 
For any $c\in\CC$, there exists a unique  structure of   $(\sigma,\rho )$-deformed $\Vc$-module  over the $h$-adically completed Yangian $\Y $ such that the module map satisfies
\beq\label{the_module_map_3}
Y_{\Y }\big(T_{[m]}^+ (u)\vac,z\big)= T_1^-(z+u_1)\ldots T_m^-(z+u_m). 
\eeq
\end{thm}

\begin{rem}
Theorem \ref{kor_tm_tm_Yang} can be reformulated as follows. By \cite[Prop. 4.2]{JKMY}, the $h$-adically completed dual Yangian $\Yd$ can be equipped with the structure of $h$-adic quantum vertex algebra $\mathcal{V}^+(\mathfrak{g}_N)$, such that its braiding, which governs \eqref{eslokaliti}, is   the map $\sigma$  given by  \eqref{sigma_formula}. Moreover, its vertex operator map is the map $Y^\rho$  given by \eqref{Yrhoo}. Thus,  $(\sigma, 1)$ is a compatible pair in $\mathcal{V}^+(\mathfrak{g}_N)$. Hence, Theorem \ref{kor_tm_tm_Yang} asserts that the $h$-adically completed Yangian $\Y $
is equipped with a unique  structure of   $(\sigma,1 )$-deformed $\mathcal{V}^+(\mathfrak{g}_N)$-module, such that the corresponding module map is again given by \eqref{the_module_map_3}.
\end{rem}

\begin{kor}\label{kor_tm_Yang}
\upshape{(a)} The coefficients of all
$$
\mathbb{T}^{-}_{[n]}  (z)\coloneqq Y_{\Y } \big(\tm (0), z \big) \in \Y [[z^{- 1}]]
$$
generate a commutative subalgebra in  $\Y $.

\noindent\upshape{(b)} The coefficients of  
$$
\mathbb{T}^{-}_{[N]} (z) = Y_{\Yggl } \big(\tN (0), z \big)\in \Yggl [[z^{-1}]]
$$
belong to the center of   $\Yggl $.
\end{kor}

\begin{rem}
Both statements of Corollary \ref{kor_tm_Yang}   are    well-known; see \cite[Chap. 1.6, 1.14]{M_book} and \cite[Prop. 5.3]{M}. In particular, the series $\mathbb{T}^{-}_{[N]} (z)$ is the {\em quantum determinant} for $\Yggl$.
 \end{rem}

Recall that the {\em extended Yangian} $\X$ is defined as the $h$-adically completed algebra over $\CC[[h]]$ generated by the elements $t_{ij}^{(r)}$, where $i,j=1,\ldots ,N$ and $r=1,2,\ldots,$ subject to the defining relation
\beq\label{DY1DY1}
R(u-v)\ts T_1^{-}(u)\ts T_2^{-}(v) =T_2^{-}(v)\ts T_1^{-}(u)\ts R(u-v),
\eeq
where $T^-(u)$ and $t_{ij}^-(u)$ are defined as in \eqref{TDY} and \eqref{tDY}, respectfully. Thus, for $\g_N=\on_N,\spn_N$, the Yangian $\Y$ is the quotient of the extended Yangian by the first relation in \eqref{DY3}, while for $\g_N=\gl_N$, the algebras $\Xlg $ and $\Ylg$ coincide.

 \begin{kor}
Let $W$ be a topologically free $\CC[[h]]$-module and $c\in\CC$. Suppose $W$ is   a $\Y$-module and denote by $T^-(u)_W $   the action of  $T^-(u )$ on $W$. There exists a unique structure of   $(\sigma,\rho )$-deformed $\Vc$-module  over   $W$ such that the module map satisfies
$$
Y_{W}\big(T_{[n]}^+ (u)\vac,z\big)= T_1^-(z+u_1)_W\ldots T_n^-(z+u_n)_W. 
$$
On the other hand, let $W$ be a  $(\sigma,\rho )$-deformed  $\Vc$-module such that  
\beq\label{uvjet_1y}
Y_W(t_{ij}^{(-1)}\vac, z)\in   \om(W,  z^{-1}W[[z^{-1}]])\quad\text{for all }i,j=1,\ldots ,N.
\eeq
There exists a unique structure of $\X$-module on $W$ satisfying
\beq\label{uvjet_2y}
T^-(z)_W=Y_W(T^+(0)\vac, z).
\eeq
\end{kor}

\begin{prf}
The first assertion of the corollary is proved by repeating the corresponding arguments from the proof of Theorem \ref{generalized_yangian_thm}. The second assertion is verified by applying the   $\sigma$-commutativity \eqref{sigma_loc_top} to $T_{13}^+(0)T_{24}^+(0)\ot w$ with $w\in W$ and using the explicit expression \eqref{sigma_formula} for the   braiding map $\sigma$. More specifically, this  shows that $Y_W(T^+(0)\vac, z)$ satisfies the $RTT$-relation for $T^-(u)$ in \eqref{DY1}. Thus, with $Y_W(T^+(0)\vac, z)$ being of the same form as  $T^- (z)$ due to \eqref{uvjet_1y}, we conclude that \eqref{uvjet_2y} defines a unique structure of $\X$-module on $W$, as required.
\end{prf}

\section{Reflection algebras  as   \texorpdfstring{$(\mu, \nu)$-deformed}{(mu,nu)-deformed} \texorpdfstring{$\mathcal{V}^{\mathrm{ crit}  }(\mathfrak{gl}_N)$}{Vcrit(glN)}-modules}\label{section_07}

We start by following \cite[Sect. 2]{MR} to define a certain class of reflection algebras. As with the generalized Yangians in the previous section, we again modify the original definition so that it is given over the ring $\CC[[h]]$. 
Consider the matrix
$$
B(u)=T(u) \ts T(-u)^{-1}\in\ndo\CC^N \ot \Yggl[[u^{-1}]].
$$
It satisfies the {\em reflection equation}
\beq\label{reflection_equation}
R(u-v)\ts B_1(u)\ts R(u+v)\ts B_2(v)=
B_2(v)\ts R(u+v)\ts  B_1(u)\ts R(u-v)
\eeq
and the {\em unitarity condition}
\beq\label{unitarity_condtion}
B(u)\ts B(-u)=1.
\eeq
The matrix $B(u)$ can be written as
$$
B(u)=\sum_{i,j=1}^N e_{ij}\ot b_{ij}(u),\quad\text{where}\quad
b_{ij}(u)=\delta_{ij}+h\sum_{r=1}^\infty b_{ij}^{(r)} u^{-r} 
$$
for some elements $b_{ij}^{(r)}\in\Yggl$.
Denote by $\BB=\mathcal{B}(N,0)$ the $h$-adically completed subalgebra of the Yangian generated by the elements $b_{ij}^{(r)}$, where $i,j=1,\ldots ,N$ and $r=1,2,\ldots .$  These generators, along with the relations in \eqref{reflection_equation} and \eqref{unitarity_condtion}, form a presentation of the algebra $\BB$; see \cite[Thm. 3.1]{MR}. We refer to $\BB $ as {\em reflection algebra}.

For any $n\geqslant 1$ and the family of variables $u=(u_1,\ldots ,u_n)$ we write
\beq\label{b-novi}
B_{[n]}  (u|z)=
\prod_{i=1,\ldots ,n}^{\longrightarrow} \left(B_i(z+u_i)\ts \R_{i\ts i+1}(2z+u_i+u_{i+1})\ldots \R_{i n}(2z+u_i+u_{n})\right).
\eeq
For example, 
for $B_i=B_i(z+u_i)$ and $\R_{ij}=\R_{ij}(2z+u_i+u_j)$,  
we have
$$
B_{[1]}(u|z)=B_1 ,\quad
B_{[2]}(u|z)=B_1 \ts \R_{12}\ts B_2 ,\quad
B_{[3]}(u|z)=B_1 \ts \R_{12}\ts \R_{13}\ts B_2(z+u_2)\ts \R_{23}\ts B_3 .
$$
By using such a notation, one can express the generalized version of   \eqref{reflection_equation}, which is a consequence of  the quantum Yang--Baxter equation \eqref{RYBE} and the reflection equation, as
\begin{align}
&R_{nm}^{12}(u|v|z_1-z_2)\ts
B_{[n]}^{13}  (u|z_1)\ts
\RR_{nm}^{12}(u|v|z_1+z_2)\ts
B_{[m]}^{23}  (v|z_2)\non\\
=&
B_{[m]}^{23}  (v|z_2)\ts
\RR_{nm}^{12}(u|v|z_1+z_2)\ts
B_{[n]}^{13}  (u|z_1)\ts
R_{nm}^{12}(u|v|z_1-z_2),\label{gen_refl_eq}
\end{align}
where the $R$-matrix product
$R_{nm}^{12}(u|v|z)$  is defined as in
\eqref{rnm12}    with the factors $\R_{ij}(z+u_i- v_{j-n})$ being replaced by $ R_{ij}(z+u_i - v_{j-n})$. 
It is worth noting that, as with the original reflection equation \eqref{reflection_equation}, the identity in \eqref{gen_refl_eq} holds if we replace all terms
$(z_1\pm z_2)$ by $( \pm z_2+z_1)$. In other words, by expanding all negative powers of $(z_1\pm z_2)$ in nonnegative powers of $z_1$ instead of $z_2$, we obtain an equivalent identity.
We shall also need the following property of  the product \eqref{b-novi}. For the variables $u=(u_1,\ldots ,u_n)$ and $v=(v_1,\ldots ,v_m)$, let $w=(u,v)=(u_1,\ldots ,u_n,v_1,\ldots ,v_m)$. Then $B_{[n+m]}  (w|z)$ decomposes as
\beq\label{b-novi-n+m}
B_{[n+m]}  (w|z)
=B^{13}_{[n]}  (u|z)\ts
\RR_{nm}^{12}(u|v|2z)\ts
B^{23}_{[m]}  (v|z).
\eeq

Consider the maps
$\gamma$, $\mu$ and $\nu$, as given by  \eqref{gamma_form}, \eqref{mu_form} and \eqref{nu_form}, respectively. 
By Proposition \ref{gen_thm_compatiblity}, the pair $(\mu,\nu)$ on $\Vccgl$ is     compatible, so one can consider the notion of $(\mu,\nu )$-deformed $\Vccgl$-module; see Definition \ref{defn_deformed_top}.
Recall that the topology $\tau$ from this definition  is required to satisfy the constraints imposed by 
 \eqref{top_submodules} and  \eqref{top_basis}. In this particular case, we define
$\BB_n =h^n\BB$ for $n\geqslant 0$ and   set $\tau$ to be the  $h$-adic topology, so that both \eqref{top_submodules} and  \eqref{top_basis} clearly hold.

\begin{thm}\label{reflection_theorem}
There exists a unique  structure of   $(\mu, \nu )$-deformed $\Vccgl$-module  over the algebra $\BB$ such that the module map is defined by
\beq\label{B-module-map}
Y_{\BB}\big(T_{[m]}^+ (u)\vac,z\big)=B_{[m]}  (u|z). 
\eeq
\end{thm}

\begin{prf}
First, let us show that the module map $Y_{\BB}(\cdot,z)$ is well-defined by \eqref{B-module-map}. It is sufficient to check that the defining relations for the dual Yangian, \eqref{DY1} with $T^+(u)$, belong to its kernel. More precisely, let $n>i$ be two positive integers,
 $u=(u_1,\ldots ,u_n)$ a family of variables and $u^\prime=(u_1,\ldots ,u_{i-1},u_{i+1},u_i,u_{i+2},\ldots,u_n)$ the $n$-tuple with $u_i$ and $u_{i+1}$ swapped.
It suffices to prove that the expression
$$
R_{i\ts i+1} \ts T_{[n]}^{+}(u)\vac
- P_{i\ts i+1}\ts T_{[n]}^{+}(u^\prime)\ts P_{i\ts i+1}\ts R_{i\ts i+1} \vac, \qquad \text{where}\qquad R_{i\ts i+1} =R_{i\ts i+1}(u_i-u_{i+1})
$$
is annihilated by $Y_{\BB}(\cdot,z)$. By using \eqref{B-module-map}, we find that it is mapped to
\beq\label{refl-well-def}
R_{i\ts i+1}\ts B_{[n]}  (u|z)
-P_{i\ts i+1}\ts B_{[n]}  (u^\prime |z)\ts P_{i\ts i+1}\ts R_{i\ts i+1}.
\eeq 
However, by examining the expression in \eqref{refl-well-def}, we see that the $R$-matrix $R_{i\ts i+1}$ in the first summand can be moved all the way to the right of $B_{[n]}  (u|z)$. Indeed, this is done by the use of  the quantum Yang--Baxter equation \eqref{RYBE} and, also, the reflection equation \eqref{reflection_equation} on the tensor factors $i$ and $i+1$. In particular,  this   turns
$$
R_{i\ts i+1}\ts B_i(z+u_i)\ts \R_{i\ts i+1}(2z+u_i+u_{i+1})\ts B_{i+1} (z+u_{i+1})
$$
to
$$
B_{i+1} (z+u_{i+1})
\ts \R_{i\ts i+1}(2z+u_i+u_{i+1})\ts
 B_i (z+u_i) \ts
R_{i\ts i+1} ,
$$
so that
$R_{i\ts i+1}  B_{[n]}  (u|z)$ becomes $P_{i\ts i+1}  B_{[n]}  (u^\prime |z)  P_{i\ts i+1}   R_{i\ts i+1}$. Hence, we conclude that \eqref{refl-well-def} is equal to zero, as required.

To finish the proof, it remains to check that the module map  $Y_{\BB}(\cdot,z)$ possesses the properties from Definition  \ref{defn_deformed_top}. The vacuum property \eqref{vac_cond_deformed} and the truncation condition \eqref{trunc_cond} clearly hold. In particular, the later is evident from the fact that the image of 
$Y_{\BB}(\cdot,z)$ belongs to $\BB[[z^{-1}]]$.
Regarding the   $\nu$-associativity, it is clear from \eqref{B-module-map} that the left-hand side of \eqref{rho_assoc_top} maps the expression $T_{[n]}^{+13}(u)  T_{[m]}^{+24}(v)(\vac\ot\vac)$ to
\beq\label{lhs-image}
B_{[n]}^{13}(u|z_2+z_0)\ts
B_{[m]}^{23}(v|z_2).
\eeq
As for the right-hand side of \eqref{rho_assoc_top},   by \eqref{Ynu-formula}, $Y^\nu(\cdot,z_0,z_2)$ maps the same expression to
 $$
    \RR_{nm}^{  12}(u|v|2z_2 +z_0-hN)^{-1}\cdotrl\left(   T_{[n]}^{+13}(u|z_0)\ts    T_{[m]}^{+23}(v)   \right) \vac .
$$
Therefore, applying the module map $Y_{\BB}(\cdot,z_2)$,   we obtain
$$
    \RR_{nm}^{  12}(u|v|2z_2 +z_0-hN)^{-1}\cdotrl B_{[n+m]}(w|z_2), 
$$
where $w=(z_0+u_1,\ldots ,z_0+u_n,v_1,\ldots ,v_m) $.
 By using \eqref{b-novi-n+m}, we rewrite this as
\begin{align}
   & \RR_{nm}^{  12}(u|v|2z_2 +z_0-hN)^{-1}\cdotrl 
		\left(
		B_{[n]}^{13}(u|z_2+z_0)\ts
		\RR_{nm}^{12}(u|v|2z_2+z_0)\ts
		B_{[m]}^{23}(v|z_2 )
		\right) \non\\
		=&\,
		B_{[n]}^{13}(u|z_2+z_0)\left(
		\RR_{nm}^{  12}(u|v|2z_2 +z_0-hN)^{-1}\cdotrl \RR_{nm}^{12}(u|v|2z_2+z_0)\right) 
		B_{[m]}^{23}(v|z_2 ).	\label{rhs-image}
\end{align}
Finally, by the crossing symmetry property \eqref{RCSYMAx}, all $R$-matrices  in \eqref{rhs-image} cancel. Thus, the resulting expression is    equal to the image \eqref{lhs-image} of the left-hand side, as required. Hence, we conclude that the   $\nu$-associativity \eqref{rho_assoc_top} holds.

Let us prove that the $\mu$-commutativity \eqref{sigma_loc_top} holds. We start by considering the image of $T_{[n]}^{+13}(u)  T_{[m]}^{+24}(v)(\vac\ot\vac)$ under the left-hand side of \eqref{sigma_loc_top}. First, we note that, by \eqref{mu_form}, its image under $\mu(-z_2+z_1,z_2)$ equals
\beq\label{mu-comm-2}
 \RR_{nm}^{  12\ts *} \cdotlr
\left(
\R_{nm}^{  12}\ts  T_{[n]}^{+13}(u)\ts
\RR_{nm}^{  12}\ts
  T_{[m]}^{+24}(v) \ts
(\R_{nm}^{  12})^{-1}  \right)(\vac\ot\vac), 
\eeq
where
$$
\R_{nm}^{  12} =\R_{nm}^{  12}(u|v|-z_2+z_1),\,
\RR_{nm}^{  12} =\RR_{nm}^{  12}(u|v|z_2+z_1 ),\,
\RR_{nm}^{  12\ts *}=\RR_{nm}^{  12}(u|v|z_2+z_1-hN)^{-1}.
$$
Next, by applying $Y_{\BB}(z_1) (1\otimes Y_{\BB}(z_2) )$ to \eqref{mu-comm-2}, we get
$$
 \RR_{nm}^{  12\ts *} \cdotlr
\left(
\R_{nm}^{  12}\ts  B_{[n]}^{13}(u|z_1)\ts
\RR_{nm}^{  12}\ts
  B_{[m]}^{23}(v|z_2) \ts
(\R_{nm}^{  12})^{-1}  \right). 
$$
By using the generalized reflection equation \eqref{gen_refl_eq}, we   rewrite the product of the first four factors inside the brackets in the above expression, thus getting
$$
 \RR_{nm}^{  12\ts *} \cdotlr
\left(
  B_{[m]}^{23}(v|z_2) \ts
	\RR_{nm}^{  12}\ts
	B_{[n]}^{13}(u|z_1)\ts
	\R_{nm}^{  12}\ts 
(\R_{nm}^{  12})^{-1}  \right)
=
 \RR_{nm}^{  12\ts *} \cdotlr
\left(
  B_{[m]}^{23}(v|z_2) \ts
	\RR_{nm}^{  12}\ts
	B_{[n]}^{13}(u|z_1) \right). 
$$
Furthermore, this simplifies to
\beq\label{mu-comm-3}
B_{[m]}^{23}(v|z_2) \left(
	\RR_{nm}^{  12\ts *} \cdotlr
	\RR_{nm}^{  12}\right)
	B_{[n]}^{13}(u|z_1) 
	=
	B_{[m]}^{23}(v|z_2) \ts
	B_{[n]}^{13}(u|z_1), 
\eeq
where the equality in \eqref{mu-comm-3} is a consequence of the crossing symmetry property \eqref{RCSYMAx}.
Finally, it is clear from \eqref{B-module-map} that the image of $T_{[n]}^{+14}(u)  T_{[m]}^{+23}(v)(\vac\ot\vac)$ under the right-hand side of \eqref{sigma_loc_top} equals the right-hand side of \eqref{mu-comm-3} , so we conclude that the $\mu$-commutativity holds, thus finishing the proof. 
\end{prf}

By combining the $\mu$-commutativity \eqref{sigma_loc_top} and the   identities   \eqref{fixed_mu} and \eqref{fixed_mu_two}, we find

\begin{kor} \label{refl_kor} 
\upshape{(a)} The coefficients of the series
$$
\mathbb{B}_{[n]}  (z)\coloneqq Y_{\BB} \big(\tm (0), z \big) \in \BB[[z^{- 1}]] \quad\text{with}\quad n=1,N 
$$
generate a commutative subalgebra in  $\BB$.

\noindent\upshape{(b)} The coefficients of  
$$
\bN (z) = Y_{\BB} \big(\bN (0), z \big)\in \BB[[z^{-1}]]
$$
belong to the center of the reflection algebra $\BB$.
\end{kor}

\begin{rem}
The commutative subalgebras  of reflection algebras in  the case of two dimensions   go   back to the pioneer work of Sklyanin \cite{S}. The series $\bN (z)$ from Corollary \ref{refl_kor} is a multiple of the well-known {\em Sklyanin determinant} $\sdet T(u)$, whose odd coefficients generate the  center of $\BB$; cf. \cite[Thm. 3.4]{MR}.
\end{rem}

At the end, we establish a connection between modules for reflection algebras and $(\mu, \nu )$-deformed $\Vccgl$-modules.
First, following  \cite[Sect. 2]{MR}, we define the {\em reflection algebra} $\BBt=\widetilde{\mathcal{B}}(N,0)$ as the $h$-adically completed algebra generated by the elements $\wtld{b}_{ij}^{(r)}$, where $i,j=1, \ldots ,N$ and $r=1,2,\ldots ,$ subject to  the {\em reflection equation}
\beq\label{reflection_equation_tilde}
R(u-v)\ts \wtld{B}_1(u)\ts R(u+v)\ts \wtld{B}_2(v)=
\wtld{B}_2(v)\ts R(u+v)\ts  \wtld{B}_1(u)\ts R(u-v),
\eeq
where
$$
\wtld{B}(u)=\sum_{i,j=1}^N e_{ij}\ot \wtld{b}_{ij}(u) \fand
\wtld{b}_{ij}(u)=\delta_{ij}+h\sum_{r=1}^\infty \wtld{b}_{ij}^{(r)} u^{-r} .
$$

\begin{kor}
Let $W$ be a topologically free $\CC[[h]]$-module. Suppose $W$ is   a $\BBt$-module and denote by $\wtld{B}(u)_W $   the action of  $\wtld{B}(u)$ on $W$. There exists a unique structure of   $(\mu, \nu )$-deformed $\Vccgl$-module  over   $W$ such that the module map satisfies
$$
Y_{W}\big(T_{[n]}^+ (u)\vac,z\big)=\wtld{B}_{[n]}  (u|z)_W, 
$$
where $\wtld{B}_{[n]}  (u|z)_W$ stands for the operator series obtained from \eqref{b-novi} by replacing all  $B_i(z+u_i)$ with $i=1,\ldots ,n$  by the  actions $\wtld{B}_i(z+u_i)_W$ on $W$. 
Conversely, let $W$ be a $(\mu, \nu )$-deformed $\Vccgl$-module such that  
\beq\label{uvjet_1}
Y_W(t_{ij}^{(-1)}\vac, z)\in   \om(W,  z^{-1}W[[z^{-1}]])\quad\text{for all }i,j=1,\ldots ,N.
\eeq
There exists a unique structure of $\BBt$-module on $W$ satisfying
\beq\label{uvjet_2}
\wtld{B}(z)_W=Y_W(T^+(0)\vac, z).
\eeq
\end{kor}

\begin{prf}
The first assertion of the corollary can be verified by repeating the arguments from the proof of Theorem \ref{reflection_theorem}. As for the second assertion, by applying the   $\mu$-commutativity \eqref{sigma_loc_top} to $T_{13}^+(0)T_{24}^+(0)\ot w$ with $w\in W$ and using the explicit expression \eqref{mu_form} for the   braiding map $\mu$, one shows that $Y_W(T^+(0)\vac, z)$ satisfies the reflection equation \eqref{reflection_equation_tilde}. Thus, with $Y_W(T^+(0)\vac, z)$ being of the same form as $\wtld{B}(z)$ due to \eqref{uvjet_1}, we conclude that \eqref{uvjet_2} defines a unique structure of $\BBt$-module on $W$.
\end{prf}

 \section*{Acknowledgment}
L. B. is member of Gruppo Nazionale per le Strutture Algebriche, Geometriche e le loro Applicazioni  (GNSAGA) of the Istituto Nazionale di Alta Matematica (INdAM) and part of the project MMNLP (Mathematical Methods in Non Linear Physics) of INFN.
L. B. was partially supported by the project Representation Theory and Applications, Bando Ateneo 2023 of Sapienza University of Rome.
S. K. is partially supported by the Croatian Science Foundation under the projects UIP-2019-04-8488 and IP-2025-02-4720  and by the project ``Implementation of cutting-edge research and its application as part of the Scientific Center of Excellence for Quantum and Complex Systems, and Representations of Lie Algebras'', Grant No. PK.1.1.10.0004, co-financed by the European Union through the European Regional Development Fund - Competitiveness and Cohesion Programme 2021--2027.
This research was supported by the European Union -- NextGenerationEU through the National Recovery and Resilience Plan 2021-2026. Institutional grant of University of Zagreb Faculty of Science Strengthening scientific production, international presence and social impact of mathematical research (IK IA 1.1.3. Impact4Math).

 \linespread{0.9}

\end{document}